\documentclass{amsart}

\usepackage[utf8]{inputenc}
\usepackage[colorlinks=true, linkcolor=black, citecolor=black, filecolor=black, urlcolor=black, pdfborder={0 0 0}, linktoc=all]{hyperref}
\usepackage{amsmath, amsthm, amssymb, mathtools, empheq, cite, enumitem}
\usepackage{todonotes, diagbox,cleveref, caption, subcaption, graphicx, pdfpages, pgf, pgfplots, tikz, pgfplots, tikz-cd,pgfmath, mathpazo, ifthen, color}
\usepackage[x11names, table]{xcolor}
\usetikzlibrary{calc, math, intersections, arrows, patterns}
\pgfplotsset{compat=1.8}
\usepackage[framemethod=TikZ]{mdframed}
\usepgfplotslibrary{groupplots}
\makeatletter

\newtheorem{theorem}{Theorem}[section]

\newtheorem{proposition}[theorem]{Proposition}
\newtheorem{lemma}[theorem]{Lemma}

\usepackage{todonotes}

\usepackage{tikz}
\usetikzlibrary{decorations.pathmorphing} 

\newcommand{\gray}{white!50!gray}
\newcommand{\black}{black}

\newcounter{pick}
\newcounter{valley}
\newcounter{valleyleft}
\newcounter{height}
\newcounter{previousheight}
\newcounter{totalheight}
\newcounter{parentheight}
\newcounter{width}
\newcounter{widthzero}
\newcounter{relativedistance}
\newcounter{toremove}
\newcounter{firstremove}

\newcommand{\nuPath}[7]{
\begin{tikzpicture}[
    scale=#3,
    pathedge/.style={line width=1.3pt}
]

\setcounter{pick}{0}
\setcounter{valley}{0}
\setcounter{height}{0}
\setcounter{previousheight}{0}
\setcounter{totalheight}{-1}

\foreach \step in {#1}{
    \addtocounter{totalheight}{1}
}

\draw[\gray,dashed] (0,0) -- (0,\thetotalheight);

\foreach \step in {#1}{

    \addtocounter{valley}{\step}

    \draw[\gray,dashed] (0,\theheight) -- (\thevalley,\theheight);

    \ifthenelse{\thepick=\thevalley}{}{
        \pgfmathparse{\thepick+1}
        \let\respick\pgfmathresult
        \foreach \x in {\respick,...,\thevalley}{
            \draw[\gray,dashed] (\x,\theheight) -- (\x,\thetotalheight);
        }
    }

    \addtocounter{pick}{\step}
    \addtocounter{height}{1}
}

\setcounter{pick}{0}
\setcounter{valley}{0}
\setcounter{height}{0}

\foreach \step in {#2}{

    \addtocounter{valley}{\step}

    \draw[pathedge]
        (\thepick,\thepreviousheight)
        --
        (\thepick,\theheight)
        --
        (\thevalley,\theheight);

    \addtocounter{pick}{\step}
    \setcounter{previousheight}{\theheight}
    \addtocounter{height}{1}
}

\foreach \a/\b/\c/\d in {#4}{
    \ifthenelse{\equal{\c}{}}{}{
        \draw[\c,line width=.5pt] (\a,\b) circle(.25);
    }
    \fill[\d,line width=.5pt] (\a,\b) circle(.15);
}

\foreach \a/\b/\c/\d/\e in {#5}{
    \draw[\d] (\a,\b) node[scale=.9,anchor=\e]{\c};
}

#7

\end{tikzpicture}
}

\newcommand{\nuTree}[7]{
\begin{tikzpicture}[
    scale=#3,
    fillnode/.style={fill=black,draw=black,circle,scale=0.8},
    circlenode/.style={draw=black,circle,scale=1},
    treeedge/.style={line width=0.9pt,black}
]

\setcounter{pick}{0}
\setcounter{valley}{0}
\setcounter{height}{0}
\setcounter{previousheight}{0}
\setcounter{totalheight}{-1}

#7
        
\foreach \step in {#1}{
    \addtocounter{totalheight}{1}
}

\draw[\gray,dashed] (0,0) -- (0,\thetotalheight);

\foreach \step in {#1}{

    \addtocounter{valley}{\step}

    \draw[\gray,dashed] (0,\theheight) -- (\thevalley,\theheight);

    \ifthenelse{\thepick=\thevalley}{}{
        \pgfmathparse{\thepick+1}
        \let\respick\pgfmathresult
        \foreach \x in {\respick,...,\thevalley}{
            \draw[\gray,dashed] (\x,\theheight) -- (\x,\thetotalheight);
        }
    }

    \addtocounter{pick}{\step}
    \addtocounter{height}{1}
}

    \setcounter{height}{0}
    \setcounter{width}{0}

    \edef\mynu{{#1}}
    \edef\mymu{{#2}}

        \pgfmathsetmacro{\stepnu}{\mynu[0]}
        \pgfmathsetmacro{\stepmu}{\mymu[0]}
        \addtocounter{width}{\stepnu}

        \foreach \x in {0,...,\stepmu}{
            \node[fillnode,scale=#3] at (\thewidth-\x,\theheight) {};
        }
        
        \draw[treeedge] (\thewidth-\stepmu,\theheight) -- (\thewidth,\theheight);

        \ifthenelse{0=\thetotalheight}{}{
        \setcounter{parentheight}{\theheight}
        \setcounter{relativedistance}{0}
        \foreach \k in {0,...,\thetotalheight}{
            \addtocounter{parentheight}{1}
            
            \pgfmathsetmacro{\substepnuabove}{\mynu[\k+1]}
            \pgfmathsetmacro{\substepmuabove}{\mymu[\k+1]}
            \addtocounter{relativedistance}{\substepnuabove-\substepmuabove}
            
            \ifthenelse{\therelativedistance>0}{}{\breakforeach} 
        }            
        \draw[treeedge] (\thewidth-\stepmu,\theheight) -- (\thewidth-\stepmu,\theparentheight);
        }

        \addtocounter{height}{1}

    \foreach \j in {1,...,\thetotalheight}{

        \pgfmathsetmacro{\stepnu}{\mynu[\j]}
        \pgfmathsetmacro{\stepmu}{\mymu[\j]}
        \addtocounter{width}{\stepnu}

        \foreach \x in {0,...,\stepmu}{
            \setcounter{relativedistance}{0}
            \setcounter{toremove}{\x}
            
            \foreach \k in {\j,...,1}{
                \pgfmathsetmacro{\substepnu}{\mynu[\k]}
                \pgfmathsetmacro{\substepmu}{\mymu[\k]}
                \pgfmathsetmacro{\substepmubelow}{\mymu[\k-1]}
                
                \ifthenelse{\k=\j}
                    {\addtocounter{relativedistance}{\x-\substepnu}}
                    {\addtocounter{relativedistance}{\substepmu-\substepnu}}

                \ifthenelse{\therelativedistance<0}
                    {\breakforeach}
                    {\addtocounter{toremove}{\substepmubelow}}
            }

            \ifthenelse{\x=0}
            {\setcounter{firstremove}{\thetoremove}}
            {}
            
            \node[fillnode,scale=#3] at (\thewidth-\thetoremove,\theheight) {};
        }

        \draw[treeedge] (\thewidth-\thetoremove,\theheight) -- (\thewidth-\thefirstremove,\theheight);
        
        \ifthenelse{\j=\thetotalheight}{}{
        \setcounter{parentheight}{\theheight}
        \setcounter{relativedistance}{0}
        \foreach \k in {\j,...,\thetotalheight}{
            \addtocounter{parentheight}{1}

            \pgfmathsetmacro{\substepnuabove}{\mynu[\k+1]}
            \pgfmathsetmacro{\substepmuabove}{\mymu[\k+1]}
            \addtocounter{relativedistance}{\substepnuabove-\substepmuabove}

            \ifthenelse{\therelativedistance>0}{}{\breakforeach} 
        }            
        \draw[treeedge] (\thewidth-\thetoremove,\theheight) -- (\thewidth-\thetoremove,\theparentheight);
        }
        
        \addtocounter{height}{1}
    }

\foreach \a/\b/\c/\d in {#4}{
    \ifthenelse{\equal{\c}{}}{}{
        \draw[\c,line width=.5pt] (\a,\b) circle(.25);
    }
    \fill[\d,line width=.5pt] (\a,\b) circle(.15);
}

\foreach \a/\b/\c/\d/\e in {#5}{
    \draw[\d] (\a,\b) node[scale=.9,anchor=\e]{\c};
}

\end{tikzpicture}
}

\newcommand{\altnuTree}[8]{
\begin{tikzpicture}[
    scale=#4,
    fillnode/.style={fill=black,draw=black,circle,scale=0.8},
    circlenode/.style={draw=black,circle,scale=1},
    treeedge/.style={line width=0.9pt,black}
]

\setcounter{valleyleft}{0}
\setcounter{height}{0}
\setcounter{previousheight}{0}
\setcounter{totalheight}{-1}

#8
        
\foreach \step in {#1}{
    \addtocounter{totalheight}{1}
}

\edef\mynu{{#1}}
\edef\mydelta{{#2}}
\edef\mymu{{#3}}

    \foreach \j in {1,...,\thetotalheight}{    
        \pgfmathsetmacro{\stepnu}{\mynu[\j]}
        \pgfmathsetmacro{\stepdelta}{\mydelta[\j-1]}
        \addtocounter{valleyleft}{\stepnu-\stepdelta}
    }
    \setcounter{valley}{\thevalleyleft}
    \pgfmathsetmacro{\stepnuzero}{\mynu[0]}
    \addtocounter{valley}{\stepnuzero}
    \setcounter{widthzero}{\thevalley}

    \draw[\gray, dashed] (\thevalleyleft,\theheight) -- (\thevalley, \theheight);
    \foreach \i in  {\thevalleyleft,...,\thevalley}{
        \draw[\gray, dashed] (\i,\theheight) -- (\i,\theheight+1);
    }
    \addtocounter{height}{1}

    \foreach \j in {1,...,\thetotalheight}{
        \pgfmathsetmacro{\stepnu}{\mynu[\j]}
        \pgfmathsetmacro{\stepdelta}{\mydelta[\j-1]}
        \addtocounter{valleyleft}{\stepdelta-\stepnu}
        \addtocounter{valley}{\stepdelta}

        \ifthenelse{\j=\thetotalheight}{
            \draw[\gray, dashed] (\thevalleyleft,\theheight) -- (\thevalley, \theheight);
            }
            {
                \draw[\gray, dashed] (\thevalleyleft,\theheight) -- (\thevalley, \theheight);
                \foreach \i in  {\thevalleyleft,...,\thevalley}{
                    \draw[\gray, dashed] (\i,\theheight) -- (\i,\theheight+1);
                }
            }
            
        \addtocounter{height}{1}
    }

    \setcounter{height}{0}
    \setcounter{width}{0}

        \pgfmathsetmacro{\stepnu}{\mynu[0]}
        \pgfmathsetmacro{\stepmu}{\mymu[0]}
        \addtocounter{width}{\thewidthzero} 

        \foreach \x in {0,...,\stepmu}{
            \node[fillnode,scale=#4] at (\thewidth-\x,\theheight) {};
        }
        
        \draw[treeedge] (\thewidth-\stepmu,\theheight) -- (\thewidth,\theheight);

        \ifthenelse{0=\thetotalheight}{}{
        \setcounter{parentheight}{\theheight}
        \setcounter{relativedistance}{0}
        \foreach \k in {0,...,\thetotalheight}{
            \addtocounter{parentheight}{1}
            
            \pgfmathsetmacro{\substepnuabove}{\mydelta[\k]} 
            \pgfmathsetmacro{\substepmuabove}{\mymu[\k+1]}
            \addtocounter{relativedistance}{\substepnuabove-\substepmuabove}
            
            \ifthenelse{\therelativedistance>0}{}{\breakforeach} 
        }            
        \draw[treeedge] (\thewidth-\stepmu,\theheight) -- (\thewidth-\stepmu,\theparentheight);
        }

        \addtocounter{height}{1}

    \foreach \j in {1,...,\thetotalheight}{

        \pgfmathsetmacro{\stepnu}{\mynu[\j]}
        \pgfmathsetmacro{\stepmu}{\mymu[\j]}
        \pgfmathsetmacro{\stepdelta}{\mydelta[\j-1]}
        \addtocounter{width}{\stepdelta} 

        \foreach \x in {0,...,\stepmu}{
            \setcounter{relativedistance}{0}
            \setcounter{toremove}{\x}
            
            \foreach \k in {\j,...,1}{
                \pgfmathsetmacro{\substepnu}{\mydelta[\k-1]} 
                \pgfmathsetmacro{\substepmu}{\mymu[\k]}
                \pgfmathsetmacro{\substepmubelow}{\mymu[\k-1]}
                
                \ifthenelse{\k=\j}
                    {\addtocounter{relativedistance}{\x-\substepnu}}
                    {\addtocounter{relativedistance}{\substepmu-\substepnu}}

                \ifthenelse{\therelativedistance<0}
                    {\breakforeach}
                    {\addtocounter{toremove}{\substepmubelow}}
            }

            \ifthenelse{\x=0}
            {\setcounter{firstremove}{\thetoremove}}
            {}
            
            \node[fillnode,scale=#4] at (\thewidth-\thetoremove,\theheight) {};
        }

        \draw[treeedge] (\thewidth-\thetoremove,\theheight) -- (\thewidth-\thefirstremove,\theheight);
        
        \ifthenelse{\j=\thetotalheight}{}{
        \setcounter{parentheight}{\theheight}
        \setcounter{relativedistance}{0}
        \foreach \k in {\j,...,\thetotalheight}{
            \addtocounter{parentheight}{1}

            \pgfmathsetmacro{\substepnuabove}{\mydelta[\k]} 
            \pgfmathsetmacro{\substepmuabove}{\mymu[\k+1]}
            \addtocounter{relativedistance}{\substepnuabove-\substepmuabove}

            \ifthenelse{\therelativedistance>0}{}{\breakforeach} 
        }            
        \draw[treeedge] (\thewidth-\thetoremove,\theheight) -- (\thewidth-\thetoremove,\theparentheight);
        }
        
        \addtocounter{height}{1}
    }

\foreach \a/\b/\c/\d in {#5}{
    \ifthenelse{\equal{\c}{}}{}{
        \draw[\c,line width=.5pt] (\a,\b) circle(.25);
    }
    \fill[\d,line width=.5pt] (\a,\b) circle(.15);
}

\foreach \a/\b/\c/\d/\e in {#6}{
    \draw[\d] (\a,\b) node[scale=.9,anchor=\e]{\c};
}

\end{tikzpicture}
}

\newcommand{\altnuShape}[8]{
\begin{tikzpicture}[
    scale=#4
]

\setcounter{valleyleft}{0}
\setcounter{height}{0}
\setcounter{previousheight}{0}
\setcounter{totalheight}{-1}

#8
        
\foreach \step in {#1}{
    \addtocounter{totalheight}{1}
}

\edef\mynu{{#1}}
\edef\mydelta{{#2}}
\edef\mymu{{#3}}

    \foreach \j in {1,...,\thetotalheight}{   
        \pgfmathsetmacro{\stepnu}{\mynu[\j]}
        \pgfmathsetmacro{\stepdelta}{\mydelta[\j-1]}
        \addtocounter{valleyleft}{\stepnu-\stepdelta}
    }
    \setcounter{valley}{\thevalleyleft}
    \pgfmathsetmacro{\stepnuzero}{\mynu[0]}
    \addtocounter{valley}{\stepnuzero}
    \setcounter{widthzero}{\thevalley}

    \draw[\black, dashed] (\thevalleyleft,\theheight) -- (\thevalley, \theheight);
    \foreach \i in  {\thevalleyleft,...,\thevalley}{
        \draw[\black, dashed] (\i,\theheight) -- (\i,\theheight+1);
    }
    \addtocounter{height}{1}

    \foreach \j in {1,...,\thetotalheight}{
        \pgfmathsetmacro{\stepnu}{\mynu[\j]}
        \pgfmathsetmacro{\stepdelta}{\mydelta[\j-1]}
        \addtocounter{valleyleft}{\stepdelta-\stepnu}
        \addtocounter{valley}{\stepdelta}

        \ifthenelse{\j=\thetotalheight}{
            \draw[\black, dashed] (\thevalleyleft,\theheight) -- (\thevalley, \theheight);
            }
            {
                \draw[\black, dashed] (\thevalleyleft,\theheight) -- (\thevalley, \theheight);
                \foreach \i in  {\thevalleyleft,...,\thevalley}{
                    \draw[\black, dashed] (\i,\theheight) -- (\i,\theheight+1);
                }
            }
            
        \addtocounter{height}{1}
    }

\foreach \a/\b/\c/\d in {#5}{
    \ifthenelse{\equal{\c}{}}{}{
        \draw[\c,line width=.5pt] (\a,\b) circle(.25);
    }
    \fill[\d,line width=.5pt] (\a,\b) circle(.15);
}

\foreach \a/\b/\c/\d/\e in {#6}{
    \draw[\d] (\a,\b) node[scale=.9,anchor=\e]{\c};
}

\end{tikzpicture}
}

\newcommand{\altnuShapegray}[8]{
\begin{tikzpicture}[
    scale=#4
]

\setcounter{valleyleft}{0}
\setcounter{height}{0}
\setcounter{previousheight}{0}
\setcounter{totalheight}{-1}

#8
        
\foreach \step in {#1}{
    \addtocounter{totalheight}{1}
}

\edef\mynu{{#1}}
\edef\mydelta{{#2}}
\edef\mymu{{#3}}

    \foreach \j in {1,...,\thetotalheight}{   
        \pgfmathsetmacro{\stepnu}{\mynu[\j]}
        \pgfmathsetmacro{\stepdelta}{\mydelta[\j-1]}
        \addtocounter{valleyleft}{\stepnu-\stepdelta}
    }
    \setcounter{valley}{\thevalleyleft}
    \pgfmathsetmacro{\stepnuzero}{\mynu[0]}
    \addtocounter{valley}{\stepnuzero}
    \setcounter{widthzero}{\thevalley}

    \draw[\gray, dashed] (\thevalleyleft,\theheight) -- (\thevalley, \theheight);
    \foreach \i in  {\thevalleyleft,...,\thevalley}{
        \draw[\gray, dashed] (\i,\theheight) -- (\i,\theheight+1);
    }
    \addtocounter{height}{1}

    \foreach \j in {1,...,\thetotalheight}{
        \pgfmathsetmacro{\stepnu}{\mynu[\j]}
        \pgfmathsetmacro{\stepdelta}{\mydelta[\j-1]}
        \addtocounter{valleyleft}{\stepdelta-\stepnu}
        \addtocounter{valley}{\stepdelta}

        \ifthenelse{\j=\thetotalheight}{
            \draw[\gray, dashed] (\thevalleyleft,\theheight) -- (\thevalley, \theheight);
            }
            {
                \draw[\gray, dashed] (\thevalleyleft,\theheight) -- (\thevalley, \theheight);
                \foreach \i in  {\thevalleyleft,...,\thevalley}{
                    \draw[\gray, dashed] (\i,\theheight) -- (\i,\theheight+1);
                }
            }
            
        \addtocounter{height}{1}
    }

\foreach \a/\b/\c/\d in {#5}{
    \ifthenelse{\equal{\c}{}}{}{
        \draw[\c,line width=.5pt] (\a,\b) circle(.25);
    }
    \fill[\d,line width=.5pt] (\a,\b) circle(.15);
}

\foreach \a/\b/\c/\d/\e in {#6}{
    \draw[\d] (\a,\b) node[scale=.9,anchor=\e]{\c};
}

\end{tikzpicture}
}

\newcommand{\altnuShapegraybonus}[9]{
\begin{tikzpicture}[
    scale=#4
]

\setcounter{valleyleft}{0}
\setcounter{height}{0}
\setcounter{previousheight}{0}
\setcounter{totalheight}{-1}

#8
        
\foreach \step in {#1}{
    \addtocounter{totalheight}{1}
}

\edef\mynu{{#1}}
\edef\mydelta{{#2}}
\edef\mymu{{#3}}

    \foreach \j in {1,...,\thetotalheight}{   
        \pgfmathsetmacro{\stepnu}{\mynu[\j]}
        \pgfmathsetmacro{\stepdelta}{\mydelta[\j-1]}
        \addtocounter{valleyleft}{\stepnu-\stepdelta}
    }
    \setcounter{valley}{\thevalleyleft}
    \pgfmathsetmacro{\stepnuzero}{\mynu[0]}
    \addtocounter{valley}{\stepnuzero}
    \setcounter{widthzero}{\thevalley}

    \draw[\gray, dashed] (\thevalleyleft,\theheight) -- (\thevalley, \theheight);
    \foreach \i in  {\thevalleyleft,...,\thevalley}{
        \draw[\gray, dashed] (\i,\theheight) -- (\i,\theheight+1);
    }
    \addtocounter{height}{1}

    \foreach \j in {1,...,\thetotalheight}{
        \pgfmathsetmacro{\stepnu}{\mynu[\j]}
        \pgfmathsetmacro{\stepdelta}{\mydelta[\j-1]}
        \addtocounter{valleyleft}{\stepdelta-\stepnu}
        \addtocounter{valley}{\stepdelta}

        \ifthenelse{\j=\thetotalheight}{
            \draw[\gray, dashed] (\thevalleyleft,\theheight) -- (\thevalley, \theheight);
            }
            {
                \draw[\gray, dashed] (\thevalleyleft,\theheight) -- (\thevalley, \theheight);
                \foreach \i in  {\thevalleyleft,...,\thevalley}{
                    \draw[\gray, dashed] (\i,\theheight) -- (\i,\theheight+1);
                }
            }
            
        \addtocounter{height}{1}
    }

\foreach \a/\b/\c/\d in {#5}{
    \ifthenelse{\equal{\c}{}}{}{
        \draw[\c,line width=.5pt] (\a,\b) circle(.25);
    }
    \fill[\d,line width=.5pt] (\a,\b) circle(.15);
}

\foreach \a/\b/\c/\d/\e in {#6}{
    \draw[\d] (\a,\b) node[scale=.9,anchor=\e]{\c};
}

#9

\end{tikzpicture}
}
\usetikzlibrary{positioning,arrows.meta,patterns.meta}
\theoremstyle{definition}
\newtheorem{definition}[theorem]{Definition}
\newtheorem{example}[theorem]{Example}

\newtheorem{remark}[theorem]{Remark}

\newcommand{\NC}[1]{\operatorname{NC}({#1})}
\newcommand{\CL}[1]{\operatorname{CL}({#1})}

\newcommand{\NCs}[1]{\operatorname{NC}_{#1}}
\newcommand{\CLnu}[1]{\operatorname{CL}_{#1}}

\title[The $\nu$-core label order and noncrossing partitions]{The $\nu$-core label order and noncrossing partitions
}
\author[M.~M\"uller]{Matthias M\"uller}
\address[M.~M\"uller]{Institute of Geometry, University of Technology Graz, Austria}
\email{matthias.mueller@tugraz.at}

\newcommand{\defn}[1]{{\color{green!50!black}\textbf{\emph{#1}}}}


\begin{document}

\begin{abstract}
We prove that the poset of $s$-noncrossing partitions is isomorphic to the core label order of the $\nu$-Tamari lattice introduced by Préville-Ratelle and Viennot.
\end{abstract}

\maketitle

\section{Introduction}

The \defn{Tamari lattice} is a fundamental partial order on Catalan objects. It can be described in several equivalent ways: as a partial order on binary trees, polygon triangulations, or on Dyck paths contained in a staircase-shaped partition. The \defn{core label order} is a poset associated with a finite lattice~\cite{reading2011noncrossing,muhle2019core}. A classical result establishes that the lattice of noncrossing partitions is isomorphic to the core label order of the Tamari lattice~\cite{reading2011noncrossing,muhle2019core}. This provides a connection between two fundamental families of Catalan objects and motivates the study of core label orders for generalizations of the Tamari lattice.
Préville-Ratelle and Viennot introduce in ~\cite{preville_vTamari_2017} the \defn{$\nu$-Tamari lattice}, a natural generalization of the Tamari lattice in which Dyck paths are replaced by an arbitrary northeast path~$\nu$. For a sequence~$s=(s_1,\dots,s_n)$ of positive integers we consider~$\nu = \nu(s) = N E^{s_1} N E^{s_2} \dots N E^{s_n}$ and establish that the core label order of the $\nu$-Tamari lattice, which we call the \defn{$\nu$-core label order}, is isomorphic to the poset of $s$-noncrossing partitions \cite{ceballos2026snoncrossing}.
\begin{theorem}\label{maintheorem}
The poset of $s$-noncrossing partitions, ordered by refinement, is isomorphic to the core label order of the $\nu$-Tamari lattice for $\nu=\nu(s)$.
\end{theorem}
\begin{figure}[h]
  \centering
  \includegraphics[width=.9\textwidth]{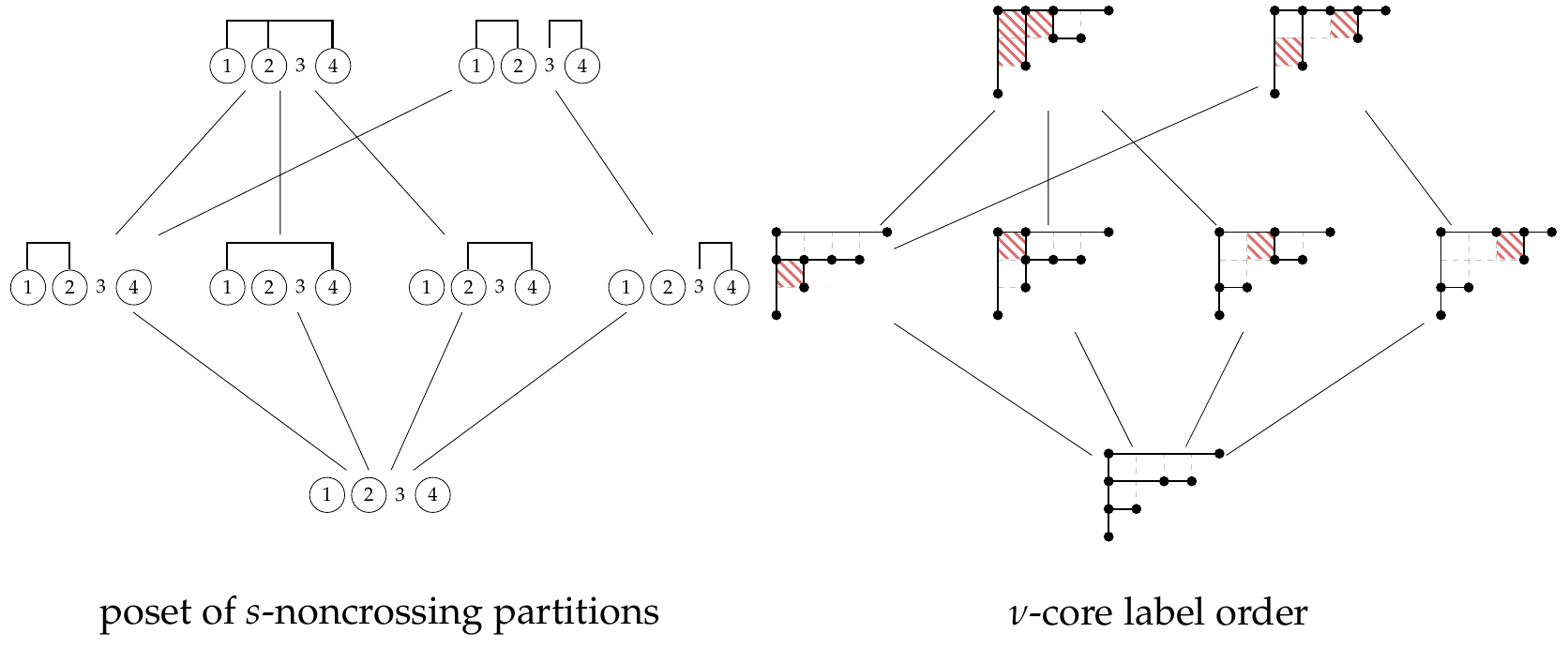}
  \caption{The poset of $s$-noncrossing partitions and the $\nu$-core label order for $s=(1,2,1)$ and $\nu=\nu(s)=NE^1NE^2 NE^1$.}
  \label{Fig_v_ncp_clo}
\end{figure}

\section{\texorpdfstring{$\nu$-trees and $s$-noncrossing partitions}{nu-trees and s-noncrossing partitions}}\label{section_bijection_trees_ncp}
In this section, we establish a bijection between $\nu$-trees and $s$-noncrossing partitions (Proposition~\ref{bijection_v_trees_ncp}). Throughout this work, $s = (s_1, \dots, s_n)$ denotes a sequence of positive integers, and $\nu = \nu(s) = NE^{s_1} \dots NE^{s_n}$.

\subsection{Background: $\nu$-trees and $s$-noncrossing partitions}
We recall the definition of the $\nu$-Tamari lattice and its elements, known as $\nu$-trees, as well as $s$-noncrossing partitions.
\subsubsection{The $\nu$-Tamari lattice}
We recall the concept of $\nu$-Tamari lattices using the conventions in~\cite{ceballos_vTamari_subword_2020}.
We denote by $\nu$ a lattice path with finitely many east and north steps. Let $F_\nu$ be the Ferrers diagram weakly above $\nu$, inside the smallest rectangle containing $\nu$. We denote by $A_\nu$ the set of lattice points weakly above $\nu$, which are inside $F_\nu$.

\begin{definition}[$\nu$-tree~\cite{ceballos_vTamari_subword_2020}]
    Two elements $ p, q \in A_\nu $ are \defn{$\nu$-incompatible}, if~$ p $ is strictly southwest (SW) of $ q $ or~$ p $ is strictly northeast (NE) of $ q $, and the smallest rectangle containing $ p $ and $q $ lies completely inside $ F_\nu $. A \defn{$\nu$-tree} is a maximal collection of pairwise~$\nu$-compatible elements in $A_\nu$. 
    Its elements are called \defn{nodes} and the top left corner is called \defn{root}.
    We associate a rooted binary tree to each~$\nu$-tree~$T$ by connecting every $p\in T$ other than the root to the next in north or west direction, see Figure~\ref{Fig_cover_relation}.
\end{definition}

\begin{definition}[$\nu$-Tamari lattice~\cite{ceballos_vTamari_subword_2020}]
Two $\nu$-trees $T$ and $T'$ are related by a right rotation, denoted by $T\lessdot T'$, if $T'$ can be obtained from $T$ by exchanging~$q \in T$ with $q' \in T'$, where~$p,r\in T,T'$, and the minimal bounding box containing $p$ and~$r$ contains no other nodes different than $q$ or $q'$.  This is illustrated in \Cref{Fig_cover_relation}. The inverse operation is called a left rotation. A node $q'$ in a $\nu$-tree $T$ is called a \defn{descent} if a left rotation can be applied to it. Similarly, it is called an \defn{ascent} if a right rotation can be applied. An example is illustrated in~\Cref{Fig_example_vTam}. The \emph{$\nu$-Tamari lattice} is the rotation poset of $\nu$-trees.
    \begin{figure}[h]
    \centering
\includegraphics[width=0.45\textwidth]{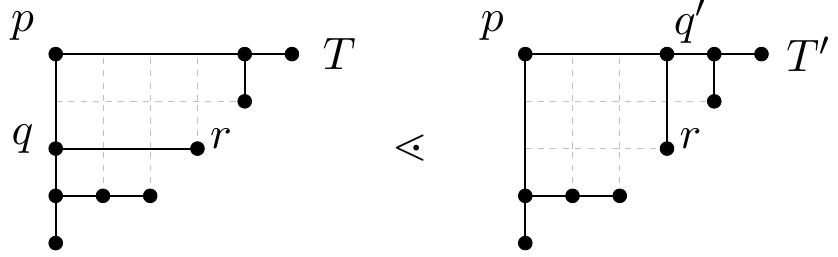} 
    \caption{right rotation of a $\nu$-tree $T$.}
    \label{Fig_cover_relation}
    \end{figure} 
    
\end{definition}


\subsubsection{$s$-noncrossing partitions}
The class of \emph{$s$-noncrossing partitions} has ben introduced in \cite{ceballos2026snoncrossing}, since it has nice enumerative properties. We recall that for a positive integer $n$, a partition $\alpha$ of $[n]:=\{1,\ldots,n\}$ is \defn{noncrossing} if there do not exist integers \(1\leq a<b<c<d\leq n\) such that \(a,c\) belong to one block of~$\alpha$ and~\(b,d\) belong to another block of~$\alpha$.
Equivalently, place $1,2,\ldots,n$ in this order on a horizontal line. For each block of $\alpha$, connect its elements by arcs drawn above the line, with consecutive elements of the block connected by an arc. The partition~$\alpha$ is noncrossing if it can be drawn such that no two such arcs cross.
We say that~$\alpha$ \defn{refines} $\alpha'$ if every block of $\alpha'$ is a union of blocks from $\alpha$; in that case, we write~$\alpha\leq \alpha'$.
The \defn{poset of noncrossing partitions} of $[n]$ with this refinement order is denoted by $\NC{n}$.

For a finite sequence of positive integers $s = (s_1, \dots, s_n)$, we let $\vert{}s\vert{} = s_1 + \dots + s_n$ and define $[s] = \{1, \dots, \vert{}s\vert{}\}$ consisting of the numbers from $1$ up to $\vert{}s\vert{}$. For each~$j \in [n]$, we define the \defn{strip}~$S_j \subseteq [s]$ as the set of indices $S_j = \{s_0 + \dots + s_{j-1} + 1, \dots, s_1 + \dots + s_j\},$ with $s_0 = 0$. The minimal element of~$S_j$ is called a \defn{first column index} of $[s]$.

\begin{definition}\label{def_s_ncp}
A subset $B \subseteq [s]$ satisfies the \defn{first column property} if
\begin{equation}\tag{FCP}\label{first_column_property}
    \text{every non-minimal element of } B \text{ is a first column index of } [s].
\end{equation}
An \defn{$s$-noncrossing partition} $\alpha$ is a noncrossing partition of $[s]$ such that each of its blocks satisfies the first column property~\eqref{first_column_property}. 
We denote by $\NCs{s}$ the \defn{poset of $s$-noncrossing partitions} ordered by refinement.
\end{definition}

\begin{example}
Let $s=(1,2,1)$, so that $[s]=\{1,2,3,4\}$, with first-column indices $1$, $2$, and $4$. The resulting poset of $s$-noncrossing partitions, $\NCs{s}$, is illustrated in~\Cref{Fig_v_ncp_clo} (left), where the first-column indices are circled. 
\end{example}

\subsection{From $\nu$-trees to $s$-noncrossing partitions} \label{section::descent_tree_cloring}
\begin{figure}[h]
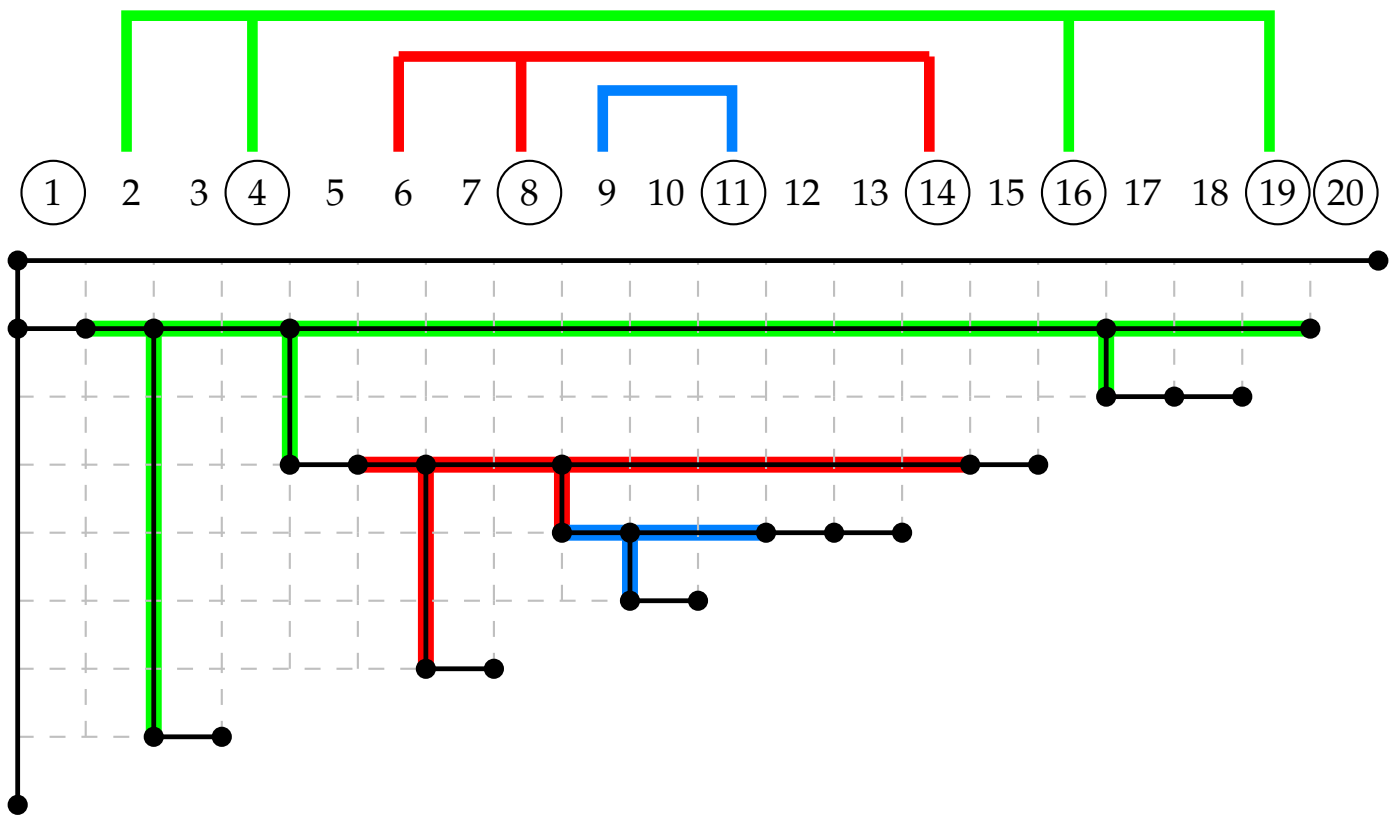

    \begin{center}
\definecolor{airforceblue}{rgb}{0.0, 0.5, 1.0}
\altnuTree{0,3,4,3,3,2,3,1,1}{3,4,3,3,2,3,1,1}{0,1,1,1,4,5,2,5,1}{0.45}{}{}{}{
\draw (1,9)  node[left, scale=0.7, black, draw, circle, minimum size=6mm, inner sep=0pt] {$1$};
\draw (2,9)  node[left, scale=0.7, black] {$2$};
\draw (3,9)  node[left, scale=0.7, black] {$3$};
\draw (4,9)  node[left, scale=0.7, black, draw, circle, minimum size=6mm, inner sep=0pt] {$4$};
\draw (5,9)  node[left, scale=0.7, black] {$5$};
\draw (6,9)  node[left, scale=0.7, black] {$6$};
\draw (7,9)  node[left, scale=0.7, black] {$7$};
\draw (8,9)  node[left, scale=0.7, black, draw, circle, minimum size=6mm, inner sep=0pt] {$8$};
\draw (9,9)  node[left, scale=0.7, black] {$9$};
\draw (10,9) node[left, scale=0.7, black] {$10$};
\draw (11,9) node[left, scale=0.7, black, draw, circle, minimum size=6mm, inner sep=0pt] {$11$};
\draw (12,9) node[left, scale=0.7, black] {$12$};
\draw (13,9) node[left, scale=0.7, black] {$13$};
\draw (14,9) node[left, scale=0.7, black, draw, circle, minimum size=6mm, inner sep=0pt] {$14$};
\draw (15,9) node[left, scale=0.7, black] {$15$};
\draw (16,9) node[left, scale=0.7, black, draw, circle, minimum size=6mm, inner sep=0pt] {$16$};
\draw (17,9) node[left, scale=0.7, black] {$17$};
\draw (18,9) node[left, scale=0.7, black] {$18$};
\draw (19,9) node[left, scale=0.7, black, draw, circle, minimum size=6mm, inner sep=0pt] {$19$};
\draw (20,9) node[left, scale=0.7, black, draw, circle, minimum size=6mm, inner sep=0pt] {$20$};

    \draw [color=green, line width=3] (1,7)--(2,7);
    \draw [color=green, line width=3] (2,7)--(2,1);
    \draw [color=green, line width=3] (2,7)--(4,7);
    \draw [color=green, line width=3] (4,5)--(4,7);
    \draw [color=green, line width=3] (4,7)--(16,7);
    \draw [color=green, line width=3] (16,6)--(16,7);
    \draw [color=green, line width=3] (16,7)--(19,7);

    \draw [color=red, line width=3] (5,5)--(14,5);
    \draw [color=red, line width=3] (6,5)--(6,2);
    \draw [color=red, line width=3] (8,5)--(8,4);

    \draw [color=airforceblue, line width=3] (8,4)--(11,4);
    \draw [color=airforceblue, line width=3] (9,4)--(9,3);
    
    \draw [color=green, line width=2] (1.6,10-.4)--(1.6,12-.4)--(18.4,12-.4)--(18.4,10-.4);
    \draw [color=green, line width=2] (3.45,10-.4)--(3.45,12-.4);
    \draw [color=green, line width=2] (15.45,10-.4)--(15.45,12-.4);

    \draw [color=red, line width=2] (5.6,10-.4)--(5.6,11);
    \draw [color=red, line width=2] (7.4,10-.4)--(7.4,11);
    \draw [color=red, line width=2] (13.4,10-.4)--(13.4,11);
    \draw [color=red, line width=2] (5.6,11)--(13.4,11);

    \draw [color=airforceblue, line width=2] (8.6,10-.4)--(8.6,10.5)--(10.5,10.5)--(10.5,10-.4);

}
\end{center}
    \caption{The descent tree coloring of a $\nu$-tree $T$ and $s$-noncrossing partition $\alpha(T)$.}
    \label{Fig_v_treecoloring}
    \end{figure}

In this subsection, we construct a map from a $\nu$-tree $T$ to an $s$-noncrossing partition $\alpha(T)$. To this end, we label the east steps of the Ferrers diagram $F_\nu$ by $1,\dots,|s|$, where $|s|$ denotes the total number of east steps of $\nu$.
Moreover, let $s_i$ is the number of east steps between the~$i$th and~$(i+1)$st north steps of $\nu$, whenever it exists. Thus, a first-column index of an $s$-noncrossing partition corresponds to an east step immediately following a north step in $\nu$. Then, the strips of $[s]$ are precisely the intervals of $\{1,\dots,|s|\}$ separated by north steps. This is illustrated in~\Cref{Fig_v_treecoloring}.

\begin{definition}\label{def_descent_tree_coloring}
The $s$-noncrossing partition $\alpha(T)$ of a $\nu$-tree $T$ is defined as follows. Let $T$ be a $\nu$-tree. Removing all nondescent nodes from $T$, different than the root, splits the graph of~$T$ into disconnected components $\widetilde{T_1},...,\widetilde{T_l}$. Let~$T_i$ be the closure of~$\widetilde{T_i}$ (include the nodes of the boundary of $\widetilde{T_i}$).
    The~$s$-noncrossing partition~$\alpha(T)$ associated to~$T$ has one block for each subtree~$T_i$, consisting of the right most column indices of its horizontal edges.
    This is illustrated in~\Cref{Fig_v_treecoloring}.        
    \end{definition}

\begin{lemma}\label{lemma_desc_tree_coloring}
    Let $T$ be a $\nu$-tree. Then, the partition $\alpha(T)$ is an $s$-noncrossing partition.
\end{lemma}
    
\begin{proof}
Let $T$ be a $\nu$-tree, and let $\alpha(T)$ be the partition obtained via the descent tree coloring of $T$. To verify Definition \ref{def_s_ncp}, we check the following:
    \begin{enumerate}
        \item The blocks $\alpha_1(T), \dots, \alpha_l(T)$  of $\alpha(T)$ are noncrossing. \label{F1}
        \item Each block $\alpha_i(T)$, $i\in[l]$ of $\alpha(T)$ satisfies the first column property \eqref{first_column_property}. \label{F2}
    \end{enumerate}
\eqref{F1}: Suppose there exist $p_i,q_i\in\alpha_i(T)$ and $p_j,q_j\in\alpha_j(T)$, $i \neq j$, such that $p_i<p_j<q_i<q_j$ or $p_j<p_i<q_j<q_i$. Without loss of generality let the horizontal edges corresponding to $p_i, q_i$ be above those for $p_j, q_j$ and label the corresponding nodes by the column index, as illustrated in \Cref{case1}. 
\begin{figure}[h]
  \centering
  \includegraphics[width=0.7\linewidth]{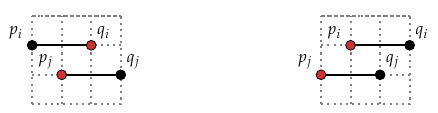}
  \caption{Illustration for the proof of Lemma~\ref{lemma_desc_tree_coloring}. Two crossing blocks induce a pair of $\nu$-incompatible nodes inside $F_\nu$.}
  \label{case1}
\end{figure}

If $p_i<p_j<q_i<q_j$, \Cref{case1} (left), then since $p_j<q_i<q_j$ the smallest rectangle containing the nodes $p_j$, $q_i$ is inside $F_\nu$. Therefore, the nodes $p_j$, $q_i$ are~$\nu$-incompatible, contradicting the definition of a $\nu$-tree.

If $p_j<p_i<q_j<q_i$, \Cref{case1} (right), then since $p_j<p_i<q_j$ the smallest rectangle containing the nodes $p_j$, $p_i$ is inside $F_\nu$. Therefore, the nodes $p_j$, $p_i$ are~$\nu$-incompatible, contradicting the definition of a $\nu$-tree.

\eqref{F2}: 
Suppose a block $\alpha_i(T)$ contains a nonfirst column index $p_2 \in \alpha_i(T)$ that is not minimal in its block, i.e., $p_2 > \min \alpha_i(T)$. Then there exists a maximal first column index $p_1 \in \alpha_i(T)$ such that $p_1 < p_2$. Since each column of a $\nu$\nobreakdash-tree contains at least one node, there exists a node associated with the right end of the first column index $\min S(p_2)$, where $S(p_2)$ denotes the strip of $F_\nu$ containing $p_2$. 

If this node is above the horizontal line between $p_1$ and $p_2$, as $q_1$ in \Cref{case3}, the nodes $p_1$ and $q_1$ are $\nu$-incompatible. On the other hand, if the node is below the horizontal line, as $q_2$ in \Cref{case3},  between $p_1$ and $p_2$, the nodes~$p_2$ and~$q_2$ are~$\nu$-incompatible, contradicting the definition of a $\nu$-tree.

\begin{figure}[h]
  \centering
  \includegraphics[width=0.2\linewidth]{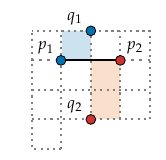}
\caption{Illustration for the proof of Lemma~\ref{lemma_desc_tree_coloring}. A block containing $p_1 < p_2$, with $p_2$ not a first-column index, induces a pair of $\nu$-incompatible nodes.}
  \label{case3}
\end{figure}

Therefore, there is a node $q$ on the horizontal line between $p_1$ and $p_2$. However, since $p_1$ is the maximal first column index smaller than $p_2$ within~$\alpha_i(T)$, it follows that there is no descent in the column of $q$. According to the definition of descent tree coloring, this implies that $p_1$ and $p_2$ cannot belong to the same block. Contradicting $p_1,p_2\in \alpha_i(T)$.
\end{proof}

The following example shows that, for a general northeast path, the descent tree coloring does not yield a bijection. In particular, to establish \Cref{maintheorem}, we require $s$ to be a sequence of positive integers.

\begin{example}
We consider $\nu = NENNE$. All $\nu$-trees and their corresponding $s$-noncrossing partitions are illustrated in \Cref{Fig_v_ncp_NENNE}. Note that distinct $\nu$-trees can yield the same $s$-noncrossing partition, which is a consequence of the consecutive north steps appearing in the lattice path $\nu$.
\begin{figure}[h]
  \centering
  \includegraphics[width=0.7\linewidth]{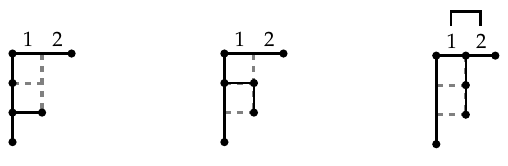}
  \caption{All $\nu$-trees and their associated $s$-noncrossing partitions via the descent tree coloring for $\nu = NENNE$.}
  \label{Fig_v_ncp_NENNE}
\end{figure}
\end{example}

\subsection{From $s$-noncrossing partitions to $\nu$-trees}
We now define the inverse construction, which sends an $s$-noncrossing partition $\alpha$ to a $\nu$-tree $T(\alpha)$. To this end, we associate a tree with each block of $\alpha$ and glue them together accordingly.


\begin{definition} \label{def_tree_by_partition}
Let $\alpha$ be an $s$-noncrossing partition with blocks $\alpha_i = \{i_1, \dots, i_k\}$ for~$i \in [l]$ such that $i_1 \le \dots \le i_k$. If $\min \alpha_i$ is a first column index, we define the tree~$T({\alpha_i})$ as the binary tree with $k-1$ internal nodes, as illustrated in \Cref{Fig_types_tree_components}~(left); otherwise, as shown in \Cref{Fig_types_tree_components}~(right). Moreover, we associate the leaves, illustrated by a white circle, with the elements $i_1, \dots, i_k$ as shown in \Cref{Fig_types_tree_components}.
\begin{figure}[h]
  \centering
  \includegraphics[width=0.7\linewidth]{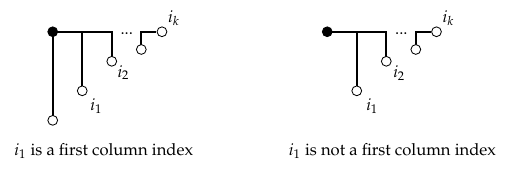}
  \caption{The two types of tree components of a $\nu$-tree for $\nu$ having no consecutive north steps.}
  \label{Fig_types_tree_components}
\end{figure}
\end{definition}

\begin{figure}[h]
  \centering
  \includegraphics[width=0.25\linewidth]{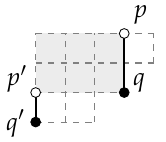}
  \caption{Illustration for the proof of Lemma~\ref{tree_components}. For $\nu$ having no consecutive north steps, no component of the the descent tree coloring is an isolated vertical edge.}
  \label{Fig_contradiction_vertical_tree_component}
\end{figure}

\begin{lemma}\label{tree_components}
For $\nu$ having no consecutive north steps, the tree components of the descent tree coloring are restricted to the two types from Definition~\ref{def_tree_by_partition}.
\end{lemma}
\begin{proof}
It suffices to show that, if $\nu$ has no consecutive north steps, no tree component (in the sense of the descent tree coloring, Definition \ref{def_descent_tree_coloring}) consists only of an isolated vertical edge.

Suppose towards a contradiction, that there exists a tree component consisting only of the edge $(p',q')$, with $p'$ above $q'$. Because in each valley of $\nu$ there is a node in $F_\nu$ and since $\nu$ has no consecutive north steps, there exists a node $q$ to the right of $p'$.    Given that the tree component $(p',q')$ is an isolated vertical edge, $q$ cannot be connected to $p'$. Consequently, $q$ must be connected to some node $p$ above it. However, this forces $p$ and $p'$ to be in $\nu$-incompatible position, as illustrated in \Cref{Fig_contradiction_vertical_tree_component}, which yields a contradiction.
\end{proof}

In the next step, we describe how to glue the trees introduced in Definition~\ref{def_tree_by_partition}.

\begin{definition}\label{def_gluing}
    Let $\alpha = \{\alpha_1, \dots, \alpha_l\}$ be an $s$-noncrossing partition, where the blocks are ordered by their minimal elements, $1 = \min \alpha_1 < \dots < \min \alpha_l$. Let $T(\alpha_1), \dots, T(\alpha_l)$ be the trees associated with the blocks of $\alpha$ by Definition~\ref{def_tree_by_partition}.    
    Let $T_1 = T(\alpha_1)$, and for each $i \ge 2$, $T_i$ is obtained by gluing the trees $T_{i-1}$ and $T(\alpha_i)$ as follows: 

\begin{itemize}
    \item If $\min(\alpha_i)$ is not a first column index, let $i'$ be the unique leaf of $T_{i-1}$ such that $i' + 1 = \min(\alpha_i)$. Then, we identify the leaf $i'$ of $T_{i-1}$ with the root of $T(\alpha_i)$. This is illustrated in \Cref{glueing} (left).
    

    \item Otherwise, if $\min(\alpha_i)$ is a first column index, $T(\alpha_i)$ is of the form \Cref{Fig_types_tree_components} (left) and contains a leaf connected to the root. We identify that leaf with the root of $T_{i-1}$. This is illustrated in \Cref{glueing} (right).
\end{itemize}

\begin{figure}[h]
  \centering
\includegraphics[width=0.45\linewidth]{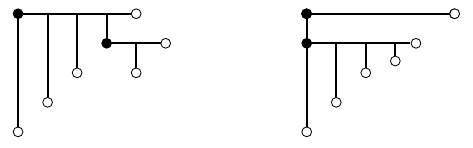}
  \caption{Illustration of the gluing step in Definition~\ref{def_gluing}. Left: $\min(\alpha_i)$ is not a first column index. Right: $\min(\alpha_i)$ is a first column index.}
  \label{glueing}
\end{figure}
\end{definition}

Having constructed the rooted binary tree $T_l$ via this recursive gluing, we now associate a $\nu(s)$-tree with it. To do this, we recall a method from \cite{ceballos_vTamari_subword_2020} to obtain a unique $\nu$-tree from a binary tree.

\begin{remark}[\textnormal{\cite[Lemma 4]{ceballos_vTamari_subword_2020}}]\label{traversal}
    We traverse the boundary of $T_l$ in counter-clockwise order from the root and associate each node $v$ of $T_l$ with the point $(x,y)$ in $A_\nu$, where $x$ is the number of~$\rightarrow$ steps before the first appearance of~$v$ in the traversal and $y$ is the number of~$\uparrow$ steps after the last appearance of $v$ in the traversal. Where the coordinate system with the root at the origin, the $x$-axis pointing to the right and the $y$-axis pointing downwards.
\end{remark}
\begin{definition}
    We define the $\nu(s)$-tree $T(\alpha)$ of an $s$-noncrossing partition $\alpha$ as the unique $\nu(s)$-tree corresponding to $T_l$ via Remark~\ref{traversal}.
\end{definition}

\begin{lemma}
    For any $s$-noncrossing partition $\alpha$, the tree $T(\alpha)$ is a $\nu$-tree, where $\nu = \nu(s)$.
\end{lemma}
\begin{proof}
We verify the following two properties:
\begin{enumerate}
    \item The number of nodes in $T(\alpha)$ matches the number of nodes in a $\nu$-tree. \label{numbernodes} 
    \item All nodes of $T(\alpha)$ form a subset of $A_\nu$ and are in $\nu$-compatible position. \label{treestructure}
\end{enumerate}
\eqref{numbernodes}: The number of nodes of a $\nu$-tree is given by $$1 + \# \text{ east steps} +\# \text{ north steps}.$$ Since $\nu = \nu(s) = NE^{s_1} \dots NE^{s_n}$, the number of nodes evaluates to $1 + |s| + n$. 
Now, for the $s$-noncrossing partition consisting entirely of singleton blocks, we have by Definition~\ref{def_tree_by_partition}, $n$ trees with $3$ nodes (\Cref{Fig_types_tree_components} (left)) and $|s|-n$ trees with $2$ nodes (\Cref{Fig_types_tree_components} (right)). Since the number of nodes decreases by $1$ in each gluing step, the total number of nodes in $T(\alpha)$ is given by $3n + 2(|s|-n) - (|s|-1) = 1 + |s| + n$.

Next, we show that the number of nodes in $T(\alpha)$ is invariant under the union of blocks in $\alpha$. Suppose we merge two blocks of $\alpha$. Then, by Definition~\ref{def_tree_by_partition}, the number of tree components decreases by one, resulting in one fewer gluing step. Since the total number of nodes in $T(\alpha)$ is given by the sum of nodes across all tree components minus the number of gluing operations, we must verify that the sum of nodes of the tree components also decreases by one. By \eqref{first_column_property}, one of two cases must occur: either we merge two blocks whose minimal elements are both first-column indices, yielding a tree component of the form shown in \Cref{Fig_types_tree_components} (left) with exactly one fewer node; or one of the blocks has a minimal element that is not a first-column index, yielding a tree component of the form shown in \Cref{Fig_types_tree_components} (right) which likewise has one fewer node. In either case, the total number of nodes in $T(\alpha)$ remains unchanged.



\eqref{treestructure}: We argue, that each node $v$ of $T(\alpha)$ is in $A_\nu$. For this it is enough to observe, that the number of $\uparrow$ steps in the counterclockwise traversal of $T(\alpha)$ after the last appearance of a node $v$ is bounded by the number of first column indices greater than the column index corresponding to $v$. Uniqueness and~$\nu$-compatibility is guaranteed by \cite[Lemma 4]{ceballos_vTamari_subword_2020}.
\end{proof}

\subsection{Bijection between $\nu$-trees and $s$-noncrossing partitions} We now show that the two constructions from above are inverse to each other.

\begin{proposition}\label{bijection_v_trees_ncp}
The set of $s$-noncrossing partitions is in bijective correspondence with the set of $\nu$-trees, where $\nu = \nu(s)$. More precisely, an $s$-noncrossing partition $\alpha$ corresponds to the $\nu$-tree $T(\alpha)$, and a $\nu$-tree $T$ corresponds to the $s$-noncrossing partition $\alpha(T)$.
\end{proposition}

\begin{proof}
This is a direct consequence of the fact that the tree components obtained via the descent tree coloring satisfy Lemma~\ref{tree_components}. Furthermore, because the gluing procedure does not introduce any new descents, the two constructions are mutually inverse.
    

\end{proof}


\section{The $\nu$-core label order}
In this section, we study the core label order of the $\nu$-Tamari lattice, called the \defn{$\nu$-core label order}, and show that it is isomorphic to the poset of $s$-noncrossing partitions (\Cref{maintheorem}).

\subsection{The core label order}
The core label order~\cite{reading2011noncrossing,muhle2019core} is a poset on the elements of a finite lattice. Originally introduced by Reading~\cite[Section~9-7.4]{Rea16}, this concept was subsequently generalized to arbitrary finite lattices by Mühle~\cite{muhle2019core,Muehle23}. We recall the definition of the core label order following\cite{Muehle23}.

Recall that a \defn{lattice} is a poset~$L$ such that any two elements~$x,y\in L$ have a \defn{meet} (i.e., greatest lower bound), denoted $x\wedge y$, and a \defn{join} (i.e., least upper bound), denoted $x\vee y$. We denote the set of cover relations by $\operatorname{Covers}(L)$ and refer to these as edges. Throughout this section, all lattices are assumed to be finite. For each element~$x \in L$, the \defn{nucleus}~is given by
$$x_\downarrow = \bigwedge_{y \in L \,:\, y \lessdot x} y,$$
and the \defn{core} is the interval~$[x_\downarrow, x]=\{u\in L \mid x_\downarrow \leq u \leq x\}$.
    A map $\lambda \colon \text{Covers}(L) \to M$ for some set $M$ is called an edge labeling of $L$. Moreover, the \defn{core label set} of $x\in L$ with respect to $\lambda$ is given by 
    \[
    \Psi_\lambda (x) = \bigl\{ \lambda(u,v) \mid x_\downarrow \leq u \lessdot v \leq x \bigr\}.
    \]
    We define a relation $\preceq_{\lambda}$ on $L$ by setting $x \preceq_{\lambda} y$ if and only if $\Psi_\lambda(x) \subseteq \Psi_\lambda(y)$ 
    and say that $\lambda$ is a \defn{core labeling} if the assignment $x \mapsto \Psi_\lambda(x)$ is injective. For any core labeling $\lambda$, the relation $\preceq_{\lambda}$ constitutes a partial order independent of $\lambda$, called the \defn{core label order} of~$L$, denoted by~$\CL{L}$.

\begin{figure}[h]
    \centering
  \scalebox{0.4}{
\begin{tikzpicture}
    \node[scale=1,circle,draw,fill=white,anchor=center,label={[font=\Huge, scale=1.2,label distance=2mm]right:$ $}] (n1) at (0,0) {};
    \draw (5,5)    node[scale=1,circle,draw,fill=white,anchor=center,label={[font=\Huge, scale=1.2,label distance=2mm]right:$ $}] (n2) {};
    \draw (8,8)    node[scale=1,circle,draw,fill=white,anchor=center,label={[font=\Huge, scale=1.2,label distance=2mm]right:$ $}] (n8) {};
\draw (-2,2)   node[circle,draw,fill=white,anchor=center,label={[font=\Huge, scale=1.2, label distance=2mm]left:$ $}] (n3) {};
    \draw (-2,5)   node[circle,draw,fill=white, anchor=center,label={[font=\Huge,scale=1.2, label distance=2mm]left:$ $}] (n4) {};
    \draw (1.5,8.5)node[scale=1,circle,draw,fill=white,anchor=center,label={[font=\Huge, scale=1.2,label distance=2mm]left: }] (n5) {};
    \draw (3+1.5,8+3+.5) node[scale=1,circle,draw,fill=white,anchor=center,label={[font=\Huge, scale=1.2,label distance=2mm]left: }] (n9) {};

    \draw [color=black, line width=1] (n5) -- node[midway, fill=white, text=red, inner sep=1pt, scale=3] {$4$} (n9);
    \draw [color=black, line width=1] (n9) -- node[midway, fill=white, text=red, inner sep=1pt, scale=3] {$1$} (n8);
    \draw [color=black, line width=1] (n8) -- node[midway, fill=white, text=red, inner sep=1pt, scale=3] {$4$} (n2);

    \draw [color=black, line width=1] (n1) -- node[midway, fill=white, text=red, inner sep=1pt, scale=3] {$3$} (n2);
    \draw [color=black, line width=1] (n2) -- node[midway, fill=white, text=red, inner sep=1pt, scale=3] {$1$} (n5);
    \draw [color=black, line width=1] (n5) -- node[midway, fill=white, text=red, inner sep=1pt, scale=3] {$3$} (n4);
    \draw [color=black, line width=1] (n4) -- node[midway, fill=white, text=red, inner sep=1pt, scale=3] {$2$} (n3);
    \draw [color=black, line width=1] (n3) -- node[midway, fill=white, text=red, inner sep=1pt, scale=3] {$1$} (n1);

\end{tikzpicture}\hspace{1cm}\begin{tikzpicture}[scale=2.7]

    \node[scale=3] (A) at (-1.5,4) {\color{black}{$\{1,2,3\}$}};
    \node[scale=3] (B) at (1,4) {\color{black}{$\{1,4\}$}};

    \node[scale=3] (a) at (-3,2) {\color{black}{$\{1\}$}};
    \node[scale=3] (b) at (-1,2) {\color{black}{$\{2\}$}};
    \node[scale=3] (c) at (1,2) {\color{black}{$\{3\}$}};
    \node[scale=3] (d) at (3,2) {\color{black}{$\{4\}$}};
    \node[scale=3] (0) at (0,0) {\color{black}{$\varnothing$}};

    \draw [color=black, line width=1] (0)--(a);
    \draw [color=black, line width=1] (0)--(b);
    \draw [color=black, line width=1] (0)--(c);
    \draw [color=black, line width=1] (0)--(d);
    
    \draw [color=black, line width=1] (a)--(A);
    \draw [color=black, line width=1] (b)--(A);
    \draw [color=black, line width=1] (c)--(A);
    \draw [color=black, line width=1] (d)--(B);
    \draw [color=black, line width=1] (a)--(B);
    
    \end{tikzpicture}
}
  \caption{Left: A finite lattice $L$ with core labeling. Right: The core label order of $L$.}
    \label{Fig_congruence_uniform_lattice}
    \end{figure}
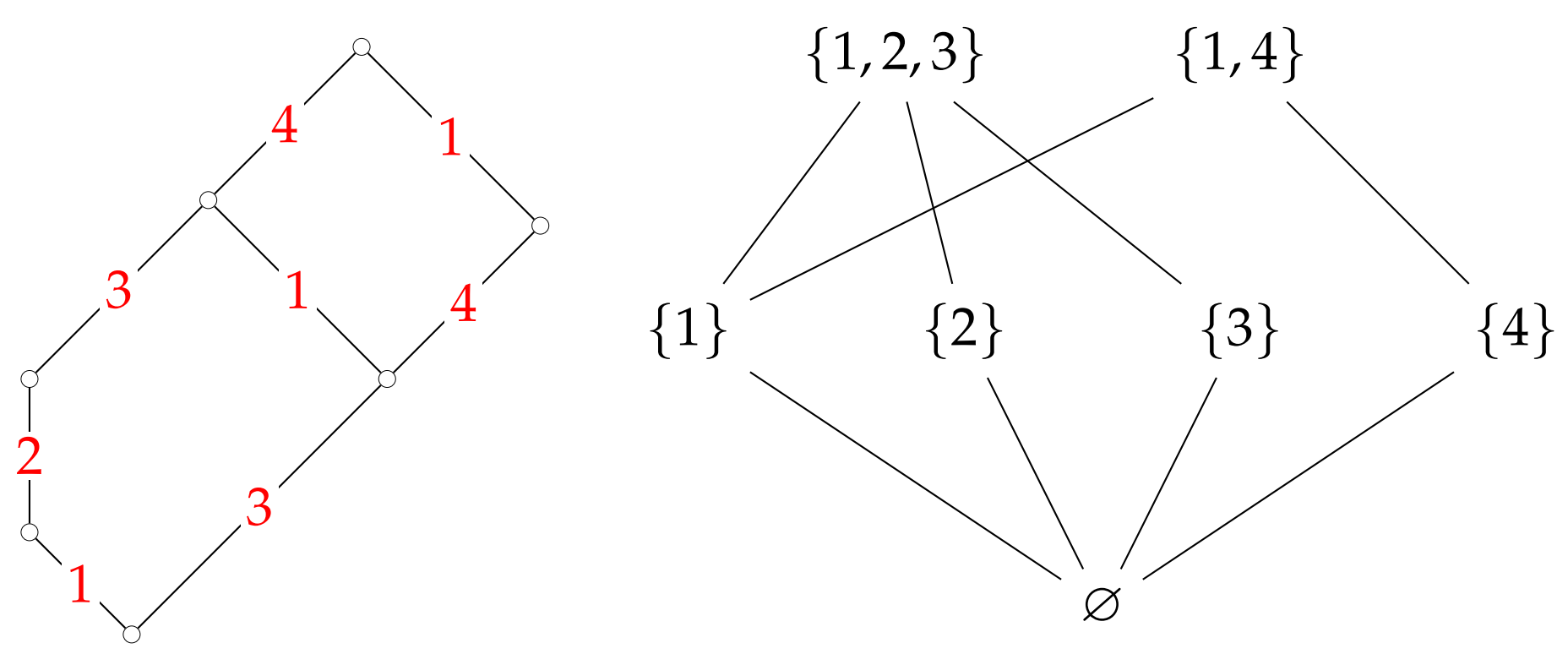

\begin{example}
Let $L$ be the lattice illustrated in~\Cref{Fig_congruence_uniform_lattice} (left), where the red labels indicate a core labeling of $L$. Identifying each element with its respective core label set and ordering these sets by containment, yields the core label order of $L$, shown in~\Cref{Fig_congruence_uniform_lattice}~(right).
\end{example}

\subsection{A core labeling of the $\nu$-Tamari lattice}
Using the description of the $\nu$-Tamari lattice in terms of $\nu$-trees, Müller introduced in~\cite{Mueller26} a labeling of the cover relations $T \lessdot T'$ by specific boxes of the Ferrers diagram $F_\nu$.

\begin{figure}[h]
    \centering
\includegraphics[width=0.6\textwidth]{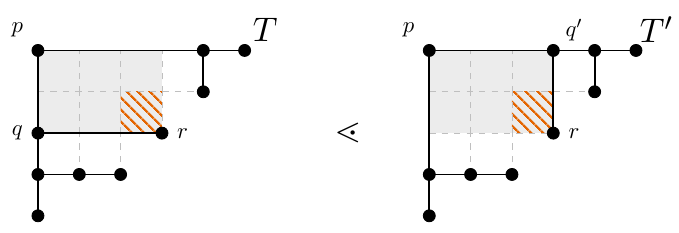} 
    \caption{A right rotation in a $\nu$-tree $T$. The box $b_{T, T'}$ labeling the cover relation $T \lessdot T'$ is striped orange.}
    \label{Fig_label}
    \end{figure}

Let $T \lessdot T'$ be a cover relation in the $\nu$-Tamari lattice described in terms of trees, and let $q$ and $q'$ be the corresponding ascent of $T$ and descent of $T'$ that are exchanged. The nodes $q$ and $q'$ form the bottom-left and top-right corners of a rectangle whose bottom-right corner $r$ is also a node of both $T$ and $T'$. We label the cover relation $T \lessdot T'$ with the box $b_{T,T'}$ of $F_\nu$ whose bottom-right corner is $r$, as illustrated in~\Cref{Fig_label}, where the rectangle is shaded gray and the box $b_{T,T'}$ is striped orange. By \cite[Lemma 2.4]{Mueller26} and \cite[Lemma 3.3]{Muehle23}, this edge labeling forms a core labeling of the $\nu$-Tamari lattice.

\subsection{The core in the $\nu$-Tamari lattice}
In this subsection, we characterize the core $[T_\downarrow, T]$ associated with a $\nu$-tree $T$. To this end, we recall the notion of bracket vectors from \cite{ceballos_vTamari_subword_2020}.

\begin{figure}[h]
    \centering
    \includegraphics[width=0.85\textwidth]{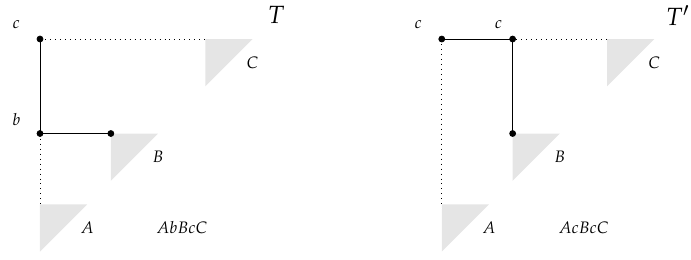} 
    \caption{A cover relation $T \lessdot T'$ between $\nu$-trees, with heights read in in-order \cite[Figure 9]{ceballos_vTamari_subword_2020}.}
    \label{Fig_v_tree_structure}
\end{figure}

\begin{figure}[h]
    \centering
    \includegraphics[width=0.8\textwidth]{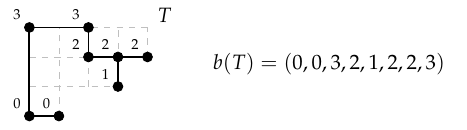} 
    \caption{A $\nu$-tree $T$ and the bracket vector $b(T)$.}
    \label{Fig8}
\end{figure}

We associate each node of a $\nu$-tree $T$ with its corresponding $y$-coordinate, from bottom to top. The~$\nu$-bracket vector $b(T)$ is then obtained by performing an \defn{in-order} traversal of $T$, also referred to as symmetric order \cite{preville_vTamari_2017}. Recursively, this traversal visits the left subtree $A$ of a root $x$, followed by $x$ itself, and concludes with the right subtree~$B$. A schematic illustration is provided in~\Cref{Fig_v_tree_structure}, which depicts a cover relation $T\lessdot T'$ together with the corresponding bracket vectors~$b(T)=AbBcC$ and $b(T')=AcBcC$. A concrete example of a~$\nu$-tree~$T$ together with its bracket vector $b(T)$ is shown in~\Cref{Fig8}. 

Ceballos, Padrol, and Sarmiento \cite{ceballos_vTamari_subword_2020} showed that the poset of $\nu$-bracket vectors under the componentwise order is isomorphic to the $\nu$-Tamari lattice. Moreover, the meet of two $\nu$-trees, $T \wedge T'$, is the unique $\nu$-tree whose bracket vector is given by the componentwise minimum of their respective bracket vectors $b(T) = (b_1, \dots, b_\ell)$ and $b(T') = (b_1', \dots, b_\ell')$ \cite[Proposition 31]{ceballos_vTamari_subword_2020}. In symbols: 
$$b(T \wedge T') = (\min \{b_1, b_1'\}, \dots, \min \{b_\ell, b_\ell'\}).$$
This provides a direct characterization of the nucleus $T_\downarrow$ of a $\nu$-tree $T$. Let $b^i = (b^i_1, \dots, b^i_\ell)$ for $i \in [k]$ denote the bracket vectors of the elements covered by $T$, where $k$ is the number of descents of $T$. Then $T_\downarrow$ is the unique $\nu$-tree such that$$b(T_\downarrow) = (\min_{1 \le i \le k} b^i_{1}, \dots, \min_{1 \le i \le k} b^i_{\ell}).$$

\begin{example}
Let $T$ be the $\nu$-tree illustrated in \Cref{ex1}. As discussed above, the bracket vector of a down cover differs from $b(T)$ in a single coordinate, transitioning from the $y$-coordinate of a descent to the $y$-coordinate of the node directly below it. Consequently, the bracket vectors $b(T)$ and $b(T^\downarrow)$ differ in exactly six positions. Specifically, the $y$-coordinates of the six descents $q_1,\dots, q_6$ in $T$ are given by $h(q_1) = h(q_2) = h(q_3) = 7$, $h(q_4) = h(q_5) = 5$, and $h(q_6) = 4$. Each entry $h(q_i)$ in $b(T)$ is replaced by a node $h(q_i')$ with $y$-coordinates $h(q_1') = 1$, $h(q_2') = 5$, $h(q_3') = 6$, $h(q_4') = 2$, $h(q_5') = 4$, and $h(q_6') = 3$. The unique $\nu$-tree $T^\downarrow$ is illustrated in \Cref{ex11}.


    \begin{figure}[ht]
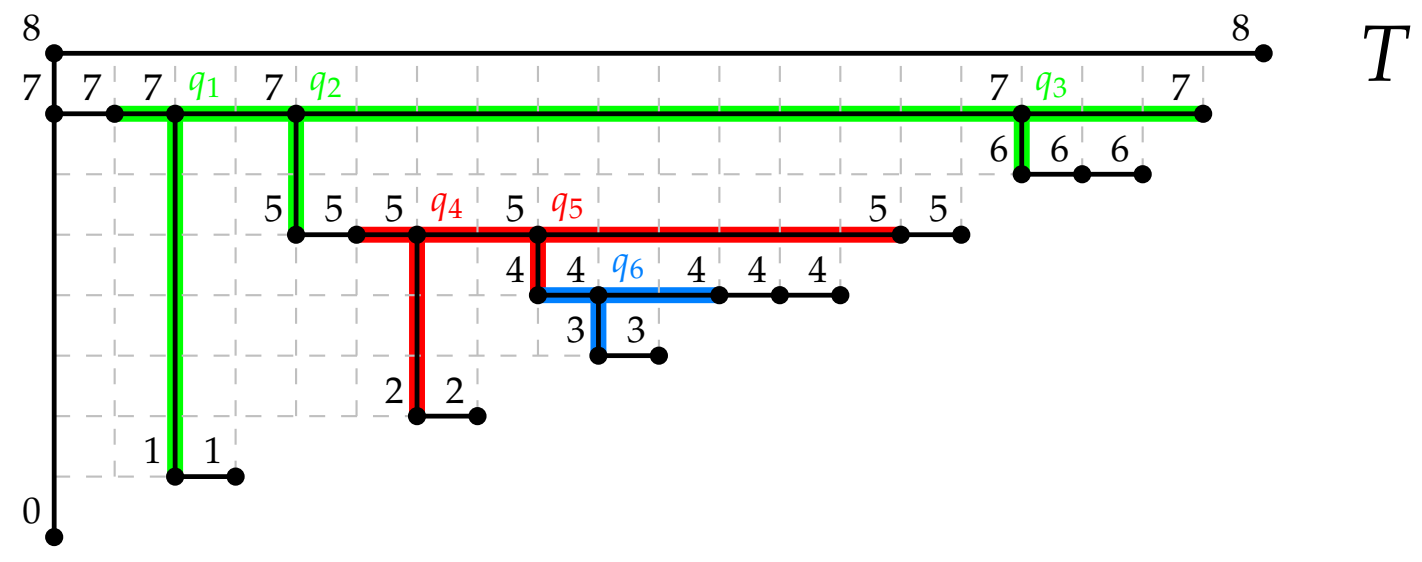

\begin{center}
\definecolor{airforceblue}{rgb}{0.0, 0.5, 1.0}
\altnuTree{0,3,4,3,3,2,3,1,1}{3,4,3,3,2,3,1,1}{0,1,1,1,4,5,2,5,1}{0.4}{}{}{}{
\draw (22,8) node[scale=1.5, black] {$T$};
\draw (0,0) node[above left, scale=0.7, black] {$0$};

\draw (2,1) node[above left, scale=0.7, black] {$1$};
\draw (3,1) node[above left, scale=0.7, black] {$1$};

\draw (6,2) node[above left, scale=0.7, black] {$2$};
\draw (7,2) node[above left, scale=0.7, black] {$2$};

\draw (9,3) node[above left, scale=0.7, black] {$3$};
\draw (10,3) node[above left, scale=0.7, black] {$3$};

\draw (8,4) node[above left, scale=0.7, black] {$4$};
\draw (9,4) node[above left, scale=0.7, black] {$4$};
\draw (11,4) node[above left, scale=0.7, black] {$4$};
\draw (12,4) node[above left, scale=0.7, black] {$4$};
\draw (13,4) node[above left, scale=0.7, black] {$4$};

\draw (4,5) node[above left, scale=0.7, black] {$5$};
\draw (5,5) node[above left, scale=0.7, black] {$5$};
\draw (6,5) node[above left, scale=0.7, black] {$5$};
\draw (8,5) node[above left, scale=0.7, black] {$5$};
\draw (14,5) node[above left, scale=0.7, black] {$5$};
\draw (15,5) node[above left, scale=0.7, black] {$5$};

\draw (16,6) node[above left, scale=0.7, black] {$6$};
\draw (17,6) node[above left, scale=0.7, black] {$6$};
\draw (18,6) node[above left, scale=0.7, black] {$6$};

\draw (0,7) node[above left, scale=0.7, black] {$7$};
\draw (1,7) node[above left, scale=0.7, black] {$7$};
\draw (2,7) node[above left, scale=0.7, black] {$7$};
\draw (4,7) node[above left, scale=0.7, black] {$7$};
\draw (5,7) node[above left, scale=0.7, green] {$q_2$};
\draw (3,7) node[above left, scale=0.7, green] {$q_1$};
\draw (17,7) node[above left, scale=0.7, green] {$q_3$};
\draw (7,5) node[above left, scale=0.7, red] {$q_4$};
\draw (9,5) node[above left, scale=0.7, red] {$q_5$};
\draw (10,4) node[above left, scale=0.7, airforceblue] {$q_6$};

\draw (16,7) node[above left, scale=0.7, black] {$7$};
\draw (19,7) node[above left, scale=0.7, black] {$7$};

\draw (0,8) node[above left, scale=0.7, black] {$8$};
\draw (20,8) node[above left, scale=0.7, black] {$8$};

    \draw [color=green, line width=3] (1,7)--(2,7);
    \draw [color=green, line width=3] (2,7)--(2,1);
    \draw [color=green, line width=3] (2,7)--(4,7);
    \draw [color=green, line width=3] (4,5)--(4,7);
    \draw [color=green, line width=3] (4,7)--(16,7);
    \draw [color=green, line width=3] (16,6)--(16,7);
    \draw [color=green, line width=3] (16,7)--(19,7);

    \draw [color=red, line width=3] (5,5)--(14,5);
    \draw [color=red, line width=3] (6,5)--(6,2);
    \draw [color=red, line width=3] (8,5)--(8,4);

    \draw [color=airforceblue, line width=3] (8,4)--(11,4);
    \draw [color=airforceblue, line width=3] (9,4)--(9,3);

}
\end{center}
\caption{A $\nu$-tree $T$ where the node labels represent their $y$-coordinates. The six descents of $T$ are labeled $q_1, \dots, q_6$. The colors indicate the tree components by the descent tree coloring.}
\label{ex1}
\end{figure}

    \begin{figure}[ht]
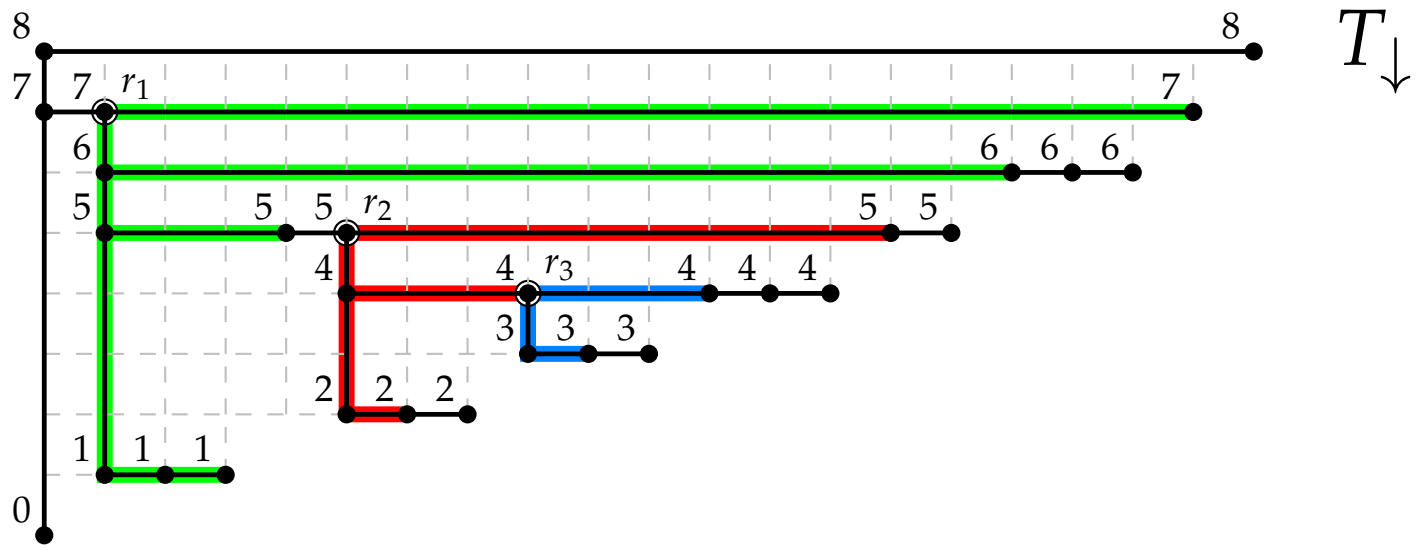

\begin{center}
\definecolor{airforceblue}{rgb}{0.0, 0.5, 1.0}
\altnuTree{0,3,4,3,3,2,3,1,1}{3,4,3,3,2,3,1,1}{0,2,2,2,4,4,3,2,1}{0.4}{}{}{}{
\draw (22,8) node[scale=1.5, black] {$T_\downarrow$};
\draw (0,0) node[above left, scale=0.7, black] {$0$};

\draw (1,1) node[above left, scale=0.7, black] {$1$};
\draw (2,1) node[above left, scale=0.7, black] {$1$};
\draw (3,1) node[above left, scale=0.7, black] {$1$};

\draw (5,2) node[above left, scale=0.7, black] {$2$};
\draw (6,2) node[above left, scale=0.7, black] {$2$};
\draw (7,2) node[above left, scale=0.7, black] {$2$};

\draw (8,3) node[above left, scale=0.7, black] {$3$};
\draw (9,3) node[above left, scale=0.7, black] {$3$};
\draw (10,3) node[above left, scale=0.7, black] {$3$};

\draw (5,4) node[above left, scale=0.7, black] {$4$};
\draw (8,4) node[above left, scale=0.7, black] {$4$};
\draw (9,4) node[above left, scale=0.7, black] {$r_3$};
\draw (11,4) node[above left, scale=0.7, black] {$4$};
\draw (12,4) node[above left, scale=0.7, black] {$4$};
\draw (13,4) node[above left, scale=0.7, black] {$4$};

\draw (1,5) node[above left, scale=0.7, black] {$5$};
\draw (4,5) node[above left, scale=0.7, black] {$5$};
\draw (5,5) node[above left, scale=0.7, black] {$5$};
\draw (6,5) node[above left, scale=0.7, black] {$r_2$};
\draw (14,5) node[above left, scale=0.7, black] {$5$};
\draw (15,5) node[above left, scale=0.7, black] {$5$};

\draw (1,6) node[above left, scale=0.7, black] {$6$};
\draw (16,6) node[above left, scale=0.7, black] {$6$};
\draw (17,6) node[above left, scale=0.7, black] {$6$};
\draw (18,6) node[above left, scale=0.7, black] {$6$};

\draw (0,7) node[above left, scale=0.7, black] {$7$};
\draw (1,7) node[above left, scale=0.7, black] {$7$};
\draw (2,7) node[above left, scale=0.7, black] {$r_1$};
\draw (19,7) node[above left, scale=0.7, black] {$7$};

\draw (0,8) node[above left, scale=0.7, black] {$8$};
\draw (20,8) node[above left, scale=0.7, black] {$8$};

    \draw [color=green, line width=3] (1,7)--(2,7);
    \draw [color=green, line width=3] (1,7)--(1,1)--(3,1);
    \draw [color=green, line width=3] (1,5)--(4,5);
    \draw [color=green, line width=3] (1,6)--(16,6);
    \draw [color=green, line width=3] (2,7)--(4,7);
    \draw [color=green, line width=3] (4,7)--(16,7);
    \draw [color=green, line width=3] (16,7)--(19,7);

    \draw [color=red, line width=3] (5,5)--(14,5);
    \draw [color=red, line width=3] (5,5)--(5,2);
    \draw [color=red, line width=3] (5,4)--(8,4);
    \draw [color=red, line width=3] (5,2)--(6,2);

    \draw [color=airforceblue, line width=3] (8,4)--(11,4);
    \draw [color=airforceblue, line width=3] (8,4)--(8,3)--(9,3);
    
    \draw (5,5) node [scale=0.5, circle, draw,fill=white,anchor=center]{};
    \draw (1,7) node [scale=0.5, circle, draw,fill=white,anchor=center]{};
    \draw (8,4) node [scale=0.5, circle, draw,fill=white,anchor=center]{};

}

\end{center}
\caption{The $\nu$-tree $T_\downarrow$ for $T$ from \Cref{ex1}. The colors indicate how the tree components of $T$ have been changed. }
\label{ex11}
\end{figure}
\end{example}

\begin{lemma}
    Fix a $\nu$-tree $T$.   
    Then the nucleus $T_\downarrow$ is the unique $\nu$-tree obtained from $T$ by applying left rotations to the descents of $T$ from left to right.
\end{lemma}
\begin{proof}
Let $q_1, \dots, q_\ell$ be the descents of $T$ ordered from left to right, which is well defined, since by $\nu$-compatibility there cannot appear two descents in the same column. It suffices to show that the bracket vector of the tree resulting from left-rotating the descents going from left to right, matches the bracket vector of the nucleus $T_\downarrow$.
Applying a left rotation to the leftmost descent $q_1$ replaces the entry of $q_1$ in the bracket vector $b(T)$ by the $y$-coordinate of the node directly below $q_1$, which is the minimum among all down covers of $T$. Proceeding to the right the same argument holds, since we choose the descents form left to right, the minimum is obtained by rotating each descent once.

\end{proof}

\begin{lemma}\label{lemma_root}
    Let $T$ be a $\nu$-tree, and let $r_1, \dots, r_\ell$ be the roots of its tree components induced by the descent tree coloring. Then, for each $\nu$-tree $\widetilde{T}' \in [T_\downarrow, T]$ such that $r_l$ is a descent for some $l \in [\ell]$, the tree $\widetilde{T}$ obtained from $\widetilde{T}'$ by a left rotation of $r_i$ (so that $\widetilde{T} \lessdot \widetilde{T} '$) is not in the core of $T$, i.e., $\widetilde{T} \notin [T_\downarrow, T]$.
\end{lemma}
\begin{proof}
Let $\widetilde{T} \lessdot \widetilde{T}'$ be a cover such that $\widetilde{T},\;\widetilde{T}' \in [T_\downarrow , T]$. Then we obtain $\widetilde{T}$ from $\widetilde{T}'$ by a left rotation, replacing $r_l$, $l\in[\ell]$, by a node $r_l'$, where the $y$-coordinate of $r_l'$ is less than that of $r_l$. Then the bracket vectors $b(\widetilde{T})$ and $b(\widetilde{T}')$ differ in exactly one entry, say the $i$th entry, where $i$ is the number of nodes read before the $y$-coordinate of~$r_l$. By Definition~\ref{def_gluing} and Proposition~\ref{bijection_v_trees_ncp}, there is a node $p$ horizontally aligned with~$r_l$ in $T_\downarrow'$. We observe that we read the $y$-coordinate of $p$ as the $i$th entry in $T_\downarrow$ and~$T$. This is because we obtain $T_\downarrow$ from $T$ by performing left rotations at the descents from left to right. Therefore, the number of entries read in-order before~$p$ does not change. Hence, $b(T_\downarrow)$ and $b(T)$ have the same $i$th entry, namely the~$y$-coordinate of $r_l$. But the $i$th entry of $\widetilde{T}'$ is $h(r_l')$. Therefore, $\widetilde{T}' \notin [T_\downarrow, T]$.
\end{proof}

\subsection{Core label sets in the $\nu$-Tamari lattice} \label{section::corelabelsets} 
In this section, we describe the core label sets in the $\nu$-Tamari lattice by a set of selected boxes inside~$F_\nu$. 

\begin{definition}
    Fix a~$\nu$-tree $T$ and let $T_1,\dots,T_\ell$ be the tree components of $T$ via the descent tree coloring, Definition~\ref{def_descent_tree_coloring}. For each tree component $T_i$ we stripe the boxes corresponding to the left rotation of its descents. Secondly, we add all boxes below an descent, where there is a marked box to the right in $T_i$. We define $\Psi'(T_i)$ as the set of striped boxes for $T_i$. We denote the set of striped boxes for $T$ by $\Psi'(T)$.

\end{definition}

\begin{example}
The set of striped boxes associated with the $\nu$-tree in \Cref{ex1} is illustrated in \Cref{ex_corelabelset}. The color of each striped box indicates the corresponding tree component. This set of boxes constitutes the core label set of the $\nu$-tree $T$, as established in the following lemma (Lemma~\ref{v_clo}).
    \begin{figure}[ht]
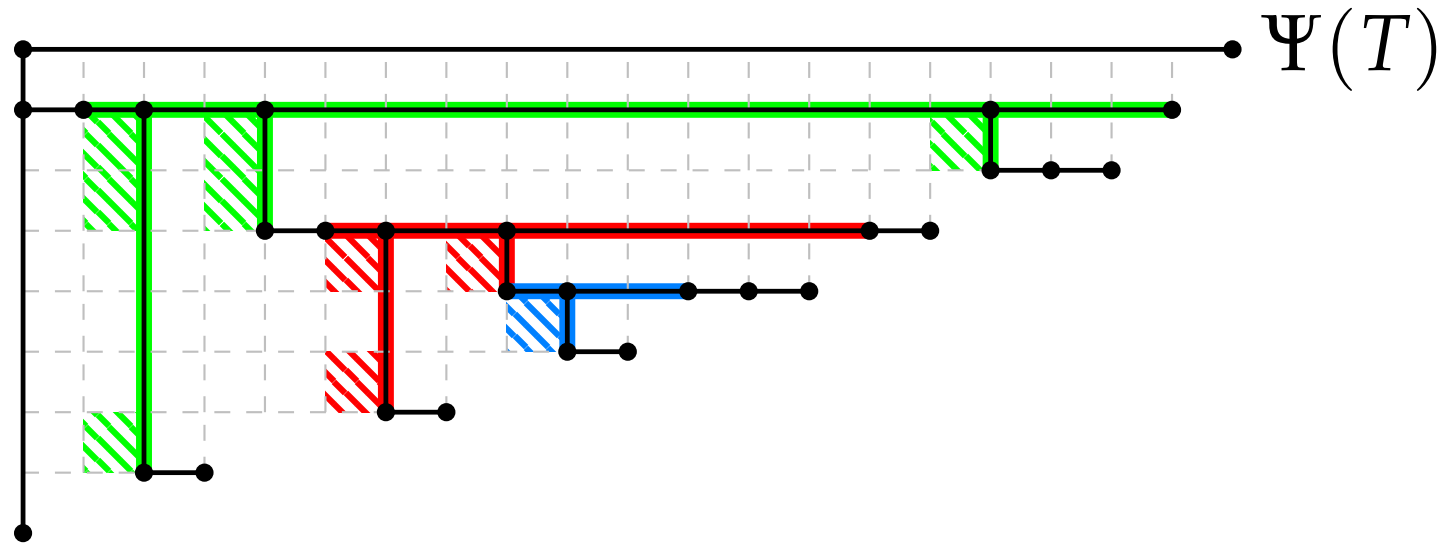

\begin{center}
\definecolor{airforceblue}{rgb}{0.0, 0.5, 1.0}

\altnuTree{0,3,4,3,3,2,3,1,1}{3,4,3,3,2,3,1,1}{0,1,1,1,4,5,2,5,1}{0.4}{}{}{}{
\draw (22,8) node[scale=1.5, black] {$\Psi (T)$};

    \draw [color=green, line width=3] (1,7)--(2,7);
    \draw [color=green, line width=3] (2,7)--(2,1);
    \draw [color=green, line width=3] (2,7)--(4,7);
    \draw [color=green, line width=3] (4,5)--(4,7);
    \draw [color=green, line width=3] (4,7)--(16,7);
    \draw [color=green, line width=3] (16,6)--(16,7);
    \draw [color=green, line width=3] (16,7)--(19,7);

    \draw [color=red, line width=3] (5,5)--(14,5);
    \draw [color=red, line width=3] (6,5)--(6,2);
    \draw [color=red, line width=3] (8,5)--(8,4);

    \draw [color=airforceblue, line width=3] (8,4)--(11,4);
    \draw [color=airforceblue, line width=3] (9,4)--(9,3);

    \fill[pattern={Lines[angle=-45, distance=3pt]}, pattern color=green, opacity=1, line width=0.4mm] (1,1) rectangle (2,2);
    \fill[pattern={Lines[angle=-45, distance=3pt]}, pattern color=green, opacity=1, line width=0.4mm] (1,5) rectangle (2,7);
    \fill[pattern={Lines[angle=-45, distance=3pt]}, pattern color=green, opacity=1, line width=0.4mm] (3,5) rectangle (4,7);
    \fill[pattern={Lines[angle=-45, distance=3pt]}, pattern color=green, opacity=1, line width=0.4mm] (15,6) rectangle (16,7);

    \fill[pattern={Lines[angle=-45, distance=3pt]}, pattern color=red, opacity=1, line width=0.4mm] (5,2) rectangle (6,3);
    \fill[pattern={Lines[angle=-45, distance=3pt]}, pattern color=red, opacity=1, line width=0.4mm] (5,4) rectangle (6,5);
    \fill[pattern={Lines[angle=-45, distance=3pt]}, pattern color=red, opacity=1, line width=0.4mm] (7,4) rectangle (8,5);

    \fill[pattern={Lines[angle=-45, distance=3pt]}, pattern color=airforceblue, opacity=1, line width=0.4mm] (8,3) rectangle (9,4);
}
\end{center}
\caption{The core label set $\Psi(T)$ for $T$ from \Cref{ex1} illustrated by striped boxes in $F_\nu$.}
\label{ex_corelabelset}
\end{figure}
\end{example}

\begin{lemma}\label{v_clo}
    For every $\nu$-tree $T$, we have $\Psi(T)=\Psi'(T)$.
    
\end{lemma}
\begin{proof}
$\Psi'(T) \subseteq \Psi(T)$: Since the poset of bracket vectors under the componentwise order is isomorphic to the $\nu$-Tamari lattice \cite{ceballos_vTamari_subword_2020}, the set of edge labels in the core~$[T_\downarrow, T]$ contains the labels of edges appearing in the longest path from $T$ to~$T_\downarrow$. This path is given by performing left rotations to the descents of$T$ from top to bottom and right to left. performing this rotations, gives precisely the boxes in~$\Psi'(T)$. Hence, all boxes in $\Psi'(T)$ are contained in the core label set $\Psi(T)$.

$\Psi'(T) \supseteq \Psi (T)$: Suppose towards a contradiction, that there is a label in the core label set $\Psi(T)$, that is not contained in $\Psi'(T)$. Then this label corresponds to a cover $\widetilde{T}\lessdot \widetilde{T}'$. Since $\widetilde{T}'$ can be obtained from $T$ by only left rotations, the desired box must lie below an descent of $T$. 
Now, $\widetilde{T}$ and $\widetilde{T}'$ contain a node in the right bottom corner of that box. Which is not a node of $T$. Hence it must come from another tree component. Where we necessarily have to rotate the root of that component. Then by Lemma~\ref{lemma_root}, $\widetilde{T} \notin [T_\downarrow ,T]$ holds. Contradicting $\widetilde{T}\lessdot \widetilde{T}'$ is a cover relation in the core $[T_\downarrow, T]$.
\end{proof}

\subsection{Proof of \texorpdfstring{\Cref{maintheorem}}{Proposition~\ref{maintheorem}}}\label{section::clo_iso}
In this section, we establish that the poset of $s$-noncrossing partitions $\NCs{s}$ is isomorphic to the $\nu$-core label order $\CLnu{\nu}$. An illustrative example of this isomorphism is provided in~\Cref{Fig_v_ncp_clo}.

\begin{figure}[h]
  \centering
\includegraphics[width=\textwidth]{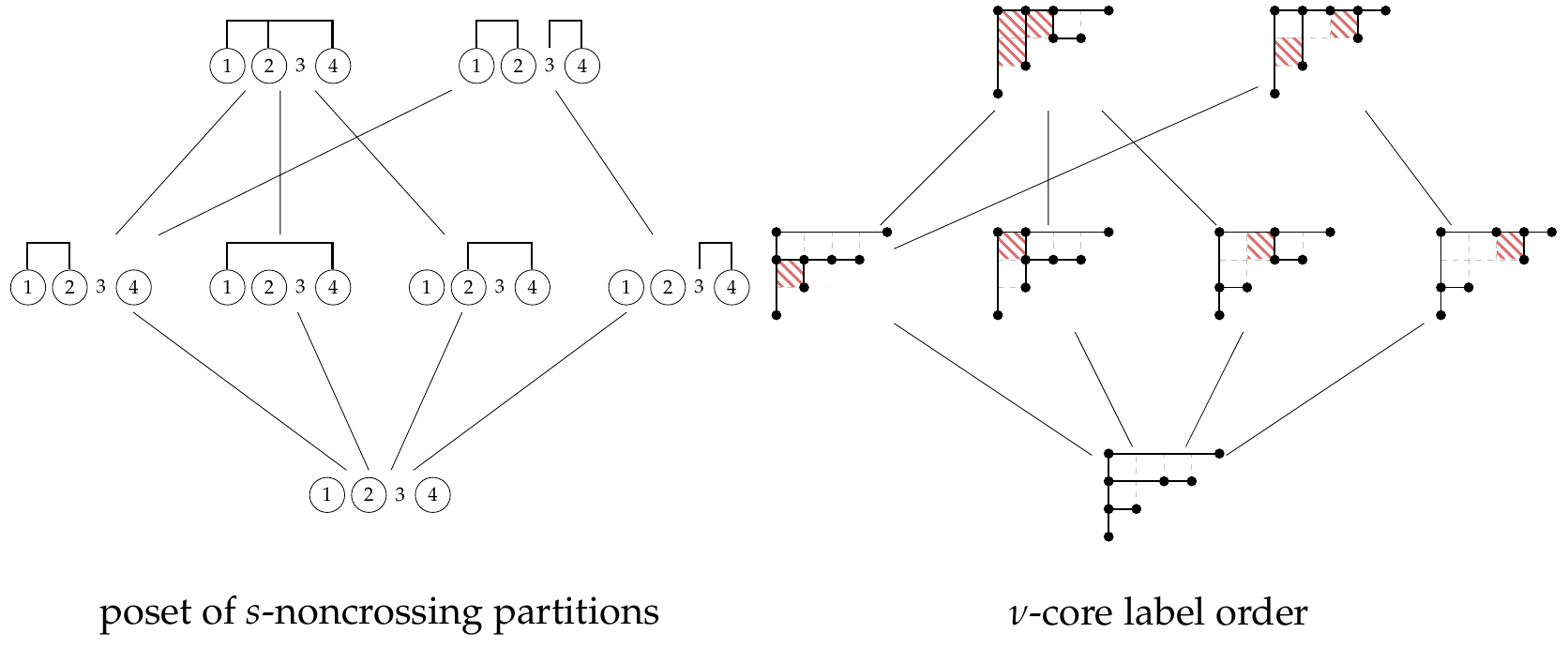}
\caption{The poset of $s$-noncrossing partitions and the $\nu$-core label order for $s=(1,2,1)$ and $\nu=\nu(s)=NENEENE$.}
    \label{Fig_v_ncp_clo}
\end{figure}

\begin{definition}
    For a lattice path $\nu$ the \defn{$\nu$-core label order} $\CLnu{\nu}$ is defined as the core label order of the $\nu$-Tamari lattice. 
\end{definition}

\begin{lemma}\label{order_preserving_1}
    For two $\nu$-trees such that $T \leq_{\CLnu{\nu}} T'$, the corresponding \(s\)-noncrossing partitions satisfy \(\alpha(T) \leq_{\NCs{s}} \alpha(T')\).
\end{lemma}
\begin{proof}
Let $T \lessdot_{\CLnu{\nu}} T'$. Suppose, for a contradiction, that $\alpha(T) \not\leq_{\NCs{s}} \alpha(T')$. Then either $\alpha(T) >_{\NCs{s}} \alpha(T')$ or the two $s$-noncrossing partitions $\alpha(T)$ and $\alpha(T')$ are incomparable. In either case, by Proposition~\ref{bijection_v_trees_ncp}, there exists at least one descent in a column of $T$ that does not appear in $T'$.Lemma~\ref{v_clo} then implies that there is an element in the core label set $\Psi(T)$ that is not contained in $\Psi(T')$. This yields $\Psi(T) \not\subseteq \Psi(T')$, which directly contradicts the initial hypothesis that $\Psi(T) \subset \Psi(T')$.
\end{proof}

The following definition characterizes $\nu$-core label sets in terms of $s$-noncrossing partitions.

\begin{definition}\label{height}

For an $s$-noncrossing partition $\alpha$ of $[s]$ with blocks  $\alpha_1, \dots, \alpha_\ell$, let the elements of each block be ordered as $\alpha_i = \{i_1, i_2, \dots, i_{n_i}\}$ such that $i_1 < i_2 < \dots < i_{n_i}$. For any block $\alpha_i$ of size $|\alpha_i| \geq 2$ and any non-maximal element $i_\ell \in \alpha_i$, the \defn{height} $h_\alpha(i_\ell)$ is defined as
$$h_\alpha(i_\ell) = |\{f \in [s] : \text{$f$ is a first column index of $[s]$ and } f \geq i_{\ell+1}\}|.$$
Let $\widetilde{\Psi}(\alpha_i)$ be the set of pairs $(i_k, h_\alpha(i_\ell))$ induced by the block $\alpha_i$:
\begin{equation}\label{boxes_ncp}
    \widetilde{\Psi}(\alpha_i) = \left\{ (i_k, h_\alpha(i_\ell)) \mid i_k, i_\ell \in \alpha_i, \text{ with } 1 \leq k \leq \ell < n_i \right\}.
\end{equation}
We define $$\widetilde{\Psi}(\alpha)=\bigcup \limits_{i=1}^n \widetilde{\Psi}(\alpha_i).$$ The set $\widetilde{\Psi}(\alpha)$ is naturally identified with a collection of boxes in $F_\nu$, where each pair~$(c, r) \in \widetilde{\Psi}(\alpha)$ represents a box located in column $c$ and row $r$ from left to right and top to bottom.
\end{definition}

\begin{lemma} \label{order_preserving_2}
    For two $s$-noncrossing partitions with $\alpha \leq_{\NCs{s}} \alpha'$, we have $\widetilde{\Psi}(\alpha) \subseteq \widetilde{\Psi}(\alpha')$.
\end{lemma}
\begin{proof}
    Let $\alpha'$ be obtain from $\alpha$ taking the union of the blocks $\{i_1,\dots, i_\ell\}$ and $\{j_1,\dots, j_k\}$, where the elements of the blocks are ordered. Since $\alpha'$ is noncrossing, we consider without loss of generality the following two cases:
    \begin{enumerate}
        \item $i_\ell \leq j_1$: The height functions of the elements $\{i_1,\dots, i_{\ell-1}\}$ and $\{j_1,\dots, j_{k-1}\}$ do not change by the Definition~\ref{height}, since the number of greater first column indices of $[s]$ does not change. Hence, $\widetilde{\Psi}(\alpha) \subsetneq \widetilde{\Psi}(\alpha')$ holds. \label{boxes_case1}
        
        \item $i_{\ell'}\leq j_1\leq \dots \leq j_k \leq i_{\ell'+1} $ for $\ell' \in \{1,\dots,\ell-1\}$: Assume without loss of generality, that $h(i_{\ell'+1})$ exists, otherwise it is the maximum and we do not associate a box with this column. Since the height of $i_{\ell'+1}$ does not change after merging, $$h_\alpha(i_{\ell'+1})=h_{\alpha'}(i_{\ell'+1}),$$ all boxes of the form $(i, h(i_{\ell'+1}))$, $i\in \{i_1,\dots,i_{\ell'}\}$ are contained in~$\widetilde{\Psi}(\alpha')$. The height of $i_{\ell'}$ is different in $\alpha$ and $\alpha'$. But the height of $i_{\ell'}$ is equal to the height of $j_k$ in $\alpha'$. This is illustrated in~\Cref{Fig_merge_boxes}. In symbols: 
        $$h_\alpha(i_{\ell'})=h_{\alpha'}(j_k).$$ 
        By Definition~\ref{height} \eqref{boxes_ncp} all the boxes in $\alpha$ at height $h_\alpha(i_{\ell'})$ are contained in~$\widetilde{\Psi}(\alpha')$ and the claim follows.
        \label{boxes_case2}
    \end{enumerate}
    \begin{figure}[h]
  \centering
\includegraphics[width=1\linewidth]{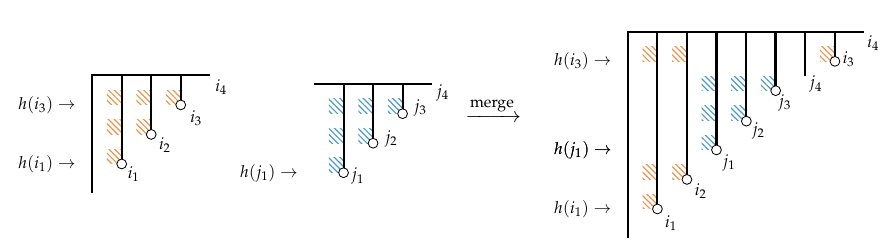}
  \caption{Merging two blocks in an $s$-noncrossing partition preserves boxes.}
  \label{Fig_merge_boxes}
\end{figure}
\end{proof}

\begin{lemma}\label{observation1}
    For an $s$-noncrossing partition $\alpha$, the height~$h_{\alpha}(i)$ for~$i \in [s]$, if it exists, is equal to the number of $\uparrow$ steps in the counter-clockwise traversal of $T(\alpha)$ that appear after the last appearance of the node $i$.
\end{lemma}
\begin{proof}
    By Definition~\ref{def_tree_by_partition}, we can identify each $\uparrow$-step in the counter-clockwise traversal of $T(\alpha)$ by a first column index greater $i$.
\end{proof}
\begin{lemma}\label{boxes_v_equality}
  Let $T$ be a $\nu$-tree, and let $\alpha = \alpha(T)$ be the $s$-noncrossing partition induced by the descent tree coloring of $T$. Then $\widetilde{\Psi}(\alpha) = \Psi(T)$ holds.
\end{lemma}
\begin{proof}
    Let $T_1,\dots,T_n$ be the tree components of $T$ by the descent tree coloring and let $\alpha(T)=\{\alpha_1, \dots, \alpha_n\}$ be the corresponding blocks in $\alpha(T)$. By the Definition~\ref{height} it is enough to show $\widetilde{\Psi}(\alpha_i) = \Psi(T_i)$. By Lemma~\ref{observation1} the~$y$-coordinate of the bottom box in column $j$ with an descent of $T_i$ coincides with the height function~$h_\alpha(j)$, hence the bottom boxes coincide. By Corollary~\ref{v_clo} all the boxes coincide.
\end{proof}

\begin{proof}[Proof of~\Cref{maintheorem}]\label{corollary_isomorphism}
By Proposition~\ref{bijection_v_trees_ncp}, $\nu$-trees are in bijection with $s$-noncrossing partitions. The equivalence~$T \leq_{\CLnu{\nu}} T' \iff \alpha(T) \leq_{\NCs{s}} \alpha(T')$ follows from Lemma~\ref{order_preserving_1} for the forward direction and Lemmas~\ref{order_preserving_2} and~\ref{boxes_v_equality} for the converse.
\end{proof}


\subsection*{Acknowledgments}
I would like to thank Cesar Ceballos for introducing me to the core label order. This work was funded by the Austrian Science Fund (FWF), project numbers  10.55776/P33278 and 10.55776/I5788.

\bibliographystyle{alpha}
\bibliography{sn-bibliography}  
\end{document}